\documentclass[10pt]{article}
\usepackage{amsmath}
\usepackage{amssymb}
\usepackage{amsfonts}
\usepackage{curves}
\numberwithin{equation}{section}
\usepackage[dvips]{graphicx}

\newtheorem{thm}{THEOREM}[section]
\newtheorem{prop}[thm]{PROPOSITION}
\newtheorem{lemma}[thm]{LEMMA}
\newtheorem{cor}[thm]{COROLLARY}

\newtheorem{remark}[thm]{REMARK}
\numberwithin{figure}{section}

\newcommand{\C}{{\mathbb C}}
\newcommand{\D}{{\mathbb D}}
\newcommand{\T}{{\mathbb T}}
\newcommand{\R}{{\mathbb R}}

\newcommand{\N}{{\mathbb N}}

\newcommand{\secpar}{{\beta}}
\newcommand{\nphip}{{\big\|(\varphi')^{t/2}\big\|}}
\newcommand\schlicht{{\mathcal S}}

\newcommand{\Cinfty}{{\mathbb C}_\infty}

\newcommand{\dz}{{\boldsymbol\partial}}
\newcommand{\bpartial}{{\boldsymbol\partial}}

\newcommand{\qed}{\vrule width 5pt height 5pt depth 0pt}

\newcommand{\calH}{{\mathcal H}}
\newcommand{\calD}{{\mathcal D}}

\newcommand{\calF}{{\mathfrak F}}

\newcommand{\calL}{{\mathcal L}}

\newcommand{\calI}{{\mathcal I}}
\newcommand{\calJ}{{\mathcal J}}
\newcommand{\calE}{{\mathcal E}}
\newcommand{\calX}{{\mathcal X}}

\newcommand{\thpar}{\sigma}
\newcommand{\fopar}{\varkappa}
\newcommand{\fipar}{K}
\newcommand{\diff}{{\mathrm d}}
\newcommand{\formderiv}{\Omega}

\newcommand{\re}{\,{\rm Re}\,}

\newcommand{\eps}{\varepsilon}
\newenvironment{proof}{\vspace{5pt}\bf Proof. \rm}{\hfill\qed\vspace{5pt}}

\newcommand{\Top}{{\mathbf T}}

\newcommand{\Pop}{{\mathbf P}}
\newcommand{\Qop}{{\mathbf Q}}

\newcommand{\Uop}{{\mathbf U}}

\newcommand{\Mop}{{\mathbf M}}

\newcommand{\spann}{\mathop{\rm span}}

\begin{document}

\font\titlefont=cmr10 scaled \magstep3
\centerline{\titlefont{Weighted Bergman spaces and the integral }}

\centerline{\titlefont{means spectrum of conformal mappings}}

\bigskip

\centerline{\sl H\aa kan Hedenmalm and Serguei Shimorin}

\centerline{\sl at the Royal Institute of Technology, Stockholm}

\bigskip

\section{Introduction}\label{se0}

\noindent{\bf The class $\schlicht$.} 
The class of univalent functions $\varphi$ from the open unit disk $\D$ into 
the complex plane $\C$, subject to the normalizations $\varphi(0)=0$ and 
$\varphi'(0)=1$, is denoted by $\schlicht$. It is classical that for 
$\varphi\in\schlicht$, we have the distortion estimates
\begin{equation}
\frac{1-|z|}{(1+|z|)^3}\le|\varphi'(z)|\le\frac{1+|z|}{(1-|z|)^3},\qquad
z\in\D.
\label{eq-koebeest}
\end{equation}
The above-mentioned estimates are sharp, as is shown by the example of a 
suitable rotation of the {\sl K\oe{}be function}
$$\kappa(z)=\frac{z}{(1-z)^2},\qquad z\in\D;$$
this function is in $\schlicht$, and maps the disk onto the plane minus 
the slit $]-\infty,-\frac14]$. After all, a simple calculation shows that
$$\kappa'(z)=\frac{1+z}{(1-z)^3},\qquad z\in\D.$$
It is of interest to better understand the sets in $\D$ where $|\varphi'(z)|$
is either large or small. For instance, $|\kappa'(z)|$ is big near the 
boundary point $z=1$, and small near $z=-1$, and elsewhere, the size is quite
modest. One way to measure the average growth or decrease is to consider 
the {\sl integral means}
$$\Mop_t[\varphi'](r)=\frac1{2\pi}\int_{-\pi}^{\pi}\big|\varphi'
\big(re^{i\theta}\big)\big|^t\,d\theta,\qquad 0<r<1,$$
where $t$ is a real parameter. It is clear from (\ref{eq-koebeest}) that
\begin{equation}
\Mop_t[\varphi'](r)=O\left(\frac1{(1-r)^\beta}\right)\,\,\,\text{as}
\,\,\,r\to1^-,
\label{eq-intbound}
\end{equation}
holds for some positive $\beta$ that depends on $t$. The {\sl infimum} of
all values of $\beta$ for which the estimate (\ref{eq-intbound}) is valid
is denoted by $\beta_\varphi(t)$. This is known as the {\sl integral means
spectral function} for $\varphi$, or simply the {\sl integral means spectrum} 
of $\varphi$. The {\sl universal integral means spectrum} for the class
$\schlicht$ is then defined by
$${\mathrm B}_\schlicht(t)=\sup_{\varphi\in\schlicht}\beta_\varphi(t).$$
Each $\beta_\varphi(t)$ is a convex function of $t$, and therefore, 
${\mathrm B}_\schlicht(t)$ is a convex function of $t$ as well.
It is a consequence of (\ref{eq-koebeest}) plus testing with $\varphi(z)=z$
that
\begin{equation}
0\le {\mathrm B}_\schlicht(t)\le \max\big\{3t,-t\big\},\qquad t\in\R.
\label{eq-trivbound}
\end{equation}
We call this the {\sl trivial bound}.

For certain $t$, the exact values of ${\mathrm B}_\schlicht(t)$ 
are known. Namely, (see \cite{FeMG})
$${\mathrm B}_\schlicht(t)=3t-1\qquad\text{for}\quad\frac25\le t<+\infty,$$
and there exists a critical value $R_{\rm CM}$, $2\le R_{\rm CM}<+\infty$ 
such that  
$${\mathrm B}_\schlicht(t)=-t-1\qquad\text{for}\quad 
-\infty<t\le-R_{\rm CM},$$ 
whereas $-t-1<{\mathrm B}_\schlicht(t)$ for $-R_{\rm CM}<t<+\infty$ 
(see \cite{CarMak}). The exact value of the universal constant $R_{\rm CM}$ 
is not known. The well-known Brennan conjecture is equivalent to the 
statement that $R_{\rm CM}=2$, which may also be expressed as 
${\mathrm B}_\schlicht(-2)=1$. 
\medskip

\noindent{\bf The class $\Sigma$.} 
We should also mention the related class $\Sigma$ of conformal maps $\varphi$ 
which map the external disk 
$$\D_e=\big\{z\in\Cinfty:\,1<|z|\le+\infty\big\}$$
into the Riemann sphere $\Cinfty=\C\cup\{\infty\}$ in such a way that
$$\varphi(z)=z+O(1),\qquad |z|\to+\infty.$$
It is classical that for $\varphi\in\Sigma$, we have the distortion estimates
\begin{equation}
\frac{|z|^2-1}{|z|^2}\le|\varphi'(z)|\le\frac{|z|^2}{|z|^2-1},\qquad
z\in\D_e.
\label{eq-koebeest'}
\end{equation}
For $\varphi\in\Sigma$, we consider the integral means
$$\Mop_t[\varphi'](r)=\frac1{2\pi}\int_{-\pi}^{\pi}\big|\varphi'
\big(re^{i\theta}\big)\big|^t\,d\theta,\qquad 1<r<+\infty,$$
and in this case, we are interested in the behavior of this quantity as
$r\to1^+$. The infimum of all $\beta$ such that
\begin{equation*}
\Mop_t[\varphi'](r)=O\left(\frac1{(r-1)^\beta}\right)\quad\text{as}
\,\,\,\,\,r\to1^+
\end{equation*}
holds is denoted by $\beta_\varphi(t)$. And ${\mathrm B}_\Sigma(t)$ -- the 
universal spectral function for the class $\Sigma$ -- is defined as the 
supremum of all $\beta_\varphi(t)$, where $\varphi$ ranges over all 
elements of $\Sigma$. 
This function ${\mathrm B}_\Sigma(t)$ is a convex function of $t$, for 
essentially the same reasons that ${\mathrm B}_\schlicht(t)$ is. The 
trivial bound of ${\mathrm B}_\Sigma(t)$ based on the pointwise estimate 
(\ref{eq-koebeest'}) is
\begin{equation}
0\le {\mathrm B}_\Sigma(t)\le |t|,\qquad t\in\R.
\label{eq-trivbound'}
\end{equation}
It is known that 
$${\mathrm B}_\Sigma(t)=|t|-1,\qquad t\in]-\infty,-R_{\rm CM}]
\cup[2,+\infty[,$$
where the constant $R_{\rm CM}$ 
is the same as before, so the 
remaining interval $[-R_{\rm CM},2]$ is what should be investigated.
\medskip

\noindent{\bf Comparison of spectra.} By analyzing the harmonic measure of
the set of points where the boundary of a simply connected set is close to 
the origin, Nikolai Makarov found in \cite{Mak1} the following relation 
between the two spectral functions:
\begin{equation}
{\mathrm B}_\schlicht(t)=\max\big\{{\mathrm B}_\Sigma(t),\,3t-1\big\},
\qquad t\in\R.
\label{eq:spectralid}
\end{equation}
We should tell the reader that Makarov's original statement deals with 
$\schlicht_b$, the class of {\sl bounded} conformal maps from $\D$ into $\C$
that preserve the origin, in place of the class $\Sigma$, but that these 
classes are sufficiently similar for the argument to carry over.

Here, we intend to study mainly the spectral function 
${\mathrm B}_\schlicht(t)$.
We shall obtain estimates that are considerably better than what has been
known up to this point. However, we have not been able to settle the 
part of the so-called Kraetzer conjecture \cite{Krae} that applies to 
${\mathrm B}_\schlicht$; this conjecture claims that
$${\mathrm B}_\Sigma(t)=\frac{t^2}4,\qquad -2\le t\le2.$$
\medskip

\noindent{\bf Bergman space methods.}
We prefer to obtain a reformulation of the definition of $\beta_\varphi(t)$
for $\varphi\in\schlicht$. It is easy to see that, for $-1<\alpha<+\infty$,
\begin{gather*}
\int_0^1 \Mop_t[\varphi'](r)\,(1-r)^{\alpha} \diff r
<+\infty\quad\Longrightarrow
\quad \Mop_t[\varphi'](r)=O\bigg(\frac1{(1-r)^{\alpha+1}}\bigg)\,\,\,\text{as}
\,\,\,r\to1^-,\\
\Mop_t[\varphi'](r)=O\bigg(\frac1{(1-r)^{\alpha+1}}\bigg)\,\,\,\text{as}
\,\,\,r\to1^-\quad\Longrightarrow\quad 
\int_0^1 \Mop_t[\varphi'](r)\,(1-r)^{\alpha+\varepsilon}\diff r<+\infty,
\end{gather*}
for each positive $\varepsilon$. For a given parameter $\alpha$ with 
$-1<\alpha<+\infty$, we now introduce the Bergman space $\calH_\alpha(\D)$, 
consisting of those holomorphic functions $f$ on $\D$ with
$$\|f\|^2_\alpha=\int_\D|f(z)|^2\,\diff A_\alpha(z)<+\infty,$$
where we use the notation
\begin{equation}
\diff A_\alpha(z)=(\alpha+1)\,\big(1-|z|^2\big)^\alpha\,\diff A(z),\qquad 
\diff A(z)=\frac{\diff x\diff y}{\pi}\,\ (z=x+iy). 
\label{eq-areael}
\end{equation}
The above expression 
defines a norm on $\calH_\alpha(\D)$ which makes it a Hilbert space.
In view of the above relationships, we have the identity
\begin{equation}
\beta_\varphi(t)=\inf\Big\{\alpha+1:\,\big(\varphi'\big)^{t/2}\in
\calH_\alpha(\D)\Big\}.
\label{eq:formula-spectrum}
\end{equation}
We think of this as a kind of ``Hilbertization'' of the problem.

In this paper, we obtain estimates of the norms 
$$\Big\|\big(\varphi'\big)^{t/2}\Big\|_\alpha$$ 
which are uniform in $\varphi\in\schlicht$; in particular, this leads to
estimates of the function ${\mathrm B}_\schlicht(t)$. 
Our methods are Bergman space techniques in combination with the classical 
tools of Geometric Function Theory, such as Gr\"onwall's area theorem. To be 
more precise, we exploit a generalization of the area theorem, due to Prawitz.
The advantage of our method is that it permits us to encode essentially 
the full strength of the area-based results, rather than just a single 
aspect thereof, such as the classical estimate ($\varphi\in\schlicht$)
\begin{equation}
\bigg|\frac{\varphi''(z)}{\varphi'(z)}-\frac{2\bar z}{1-|z|^2}\bigg|
\le\frac{4}{1-|z|^2},\qquad z\in\D,
\label{eq:Bieberbach}
\end{equation}
which is a consequence of Bieberbach's inequality
$\frac12|\varphi''(0)|=|\widehat\varphi(2)|\le2$. 
\medskip

\noindent{\bf Complex parameters in the spectral function.}
It is natural to consider the integral means spectral functions also for
complex arguments. For complex $\tau\in\C$, we define the associated 
$\tau$-integral means of $\varphi'$ by
$$\Mop_t[\varphi'](r)=\frac1{2\pi}\int_{-\pi}^{\pi}\Big|\big[\varphi'
\big(re^{i\theta}\big)\big]^{\tau}\Big|\,\diff\theta,\qquad 0<r<1,$$
for $\varphi\in\schlicht$, and by the same formula with $1<r<+\infty$
for $\varphi\in\Sigma$. The definition of the power is more delicate
this time, but we are saved by the fact that $\varphi'(z)$ is zero-free
in the disk, and we choose -- as a matter of convenience -- the branch of 
$[\varphi'(z)]^\tau$ which gives the value $1$ for $z=0$. 
This allows us to define 
$\beta_\varphi(\tau)$ just as before, and taking the suprema over the 
two classes $\schlicht$ and $\Sigma$, we obtain the universal integral
means spectral functions ${\mathrm B}_\schlicht(\tau)$ and 
${\mathrm B}_\Sigma(\tau)$ defined over $\tau\in\C$. A simple analysis of
these two functions shows that each is convex in the whole complex plane. 
Our method will supply estimates of the function 
${\mathrm B}_\schlicht(\tau)$ for complex $\tau$, but we usually do not 
stress this fact.
\medskip

\noindent{\bf Underlying ideas.} We outline the underlying philosophy of 
the paper. As we began this study of integral means spectral functions, 
we got increasingly convinced that the topic is related to the smallness of 
certain operators associated to a given conformal mapping $\varphi$. 
To get the basic idea, we suppose that 
\begin{equation*}
\sup_\varphi\int_\D \big|[\varphi'(z)]^{\tau}\big|\,\diff A_\alpha(z)<
+\infty
\label{eq:intro-1}
\end{equation*}
holds for some $\alpha$, $-1<\alpha<+\infty$, and some complex $\tau$;
the supremum runs over all $\varphi\in\schlicht$. This assumption
looks slightly stronger than the statement that ${\mathrm B}_\schlicht(\tau)
<\alpha+1$, due to the uniformity in the bound, but is most likely 
equivalent to it. We suppose that, in addition, the same estimate holds
for $-\tau$ as well:
\begin{equation*}
\sup_\varphi\int_\D \big|[\varphi'(z)]^{-\tau}\big|\,\diff A_\alpha(z)<
+\infty.
\label{eq:intro-2}
\end{equation*}
In fact, the estimate we really need is
\begin{equation*}
\sup_\varphi\,\,\,
\int_\D \big|[\varphi'(z)]^{\tau}\big|\,\diff A_\alpha(z)\times
\int_\D \big|[\varphi'(z)]^{-\tau}\big|\,\diff A_\alpha(z)<
+\infty,
\label{eq:intro-3}
\end{equation*}
which we write in the form
\begin{equation}
\sup_\varphi\,\,\,
\Big\langle\big|[\varphi'(z)]^{\tau}\big|\Big\rangle_{\D,\alpha}
\Big\langle\big|[\varphi'(z)]^{\tau}\big|^{-1}\Big\rangle_{\D,\alpha}
+\infty,
\label{eq:intro-4}
\end{equation}
where the notation
$$\big\langle f\big\rangle_{Q,\alpha}=\frac1{|Q|_\alpha}\int_Q
f(z)\,\diff A_\alpha(z)$$
is used for the $\diff A_\alpha$-average of $f$ on the subset $Q$ of $\D$; 
here, $|Q|_\alpha$ is the $\diff A_\alpha$-area of $Q$. 
We now use the fact that for each $z_0\in\D$, the function
$$z\mapsto\frac{\varphi\left(\frac{z+z_0}{1+\bar z_0\,z}\right)-
\varphi(z_0)}{(1-|z_0|^2)\,\varphi'(z_0)}$$
is an element of $\schlicht$, plug it into (\ref{eq:intro-4}) in place of
$\varphi$, and make an appropriate M\oe{}bius shift of coordinates in $\D$. 
It then follows that
\begin{equation}
\sup_{\varphi,Q}\,\,\,
\Big\langle\big|[\varphi'(z)]^{\tau}\big|\Big\rangle_{Q,\alpha}
\Big\langle\big|[\varphi'(z)]^{\tau}\big|^{-1}\Big\rangle_{Q,\alpha}
+\infty,
\label{eq:intro-5}
\end{equation}
where the supremum runs over all $\varphi\in\schlicht$ and all Carleson
``squares'' $Q$in $\D$. Condition (\ref{eq:intro-5}) is of 
$\diff A_\alpha$-area Muckenhoupt (or B\'ekoll\'e) type. In the limit case 
$\alpha=-1$, when $\diff A_\alpha$ degenerates to arc length measure on the 
unit circle $\T$, the Muckenhoupt $(A_2)$ condition on the positive weight 
$\omega$, which reads
$$\sup_Q\,\,
\big\langle\omega\big\rangle_{Q,-1}
\big\langle\omega^{-1}\big\rangle_{Q,-1}
+\infty,$$
is -- by the celebrated Helson-Szeg\"o theorem \cite{HelSzeg} --
equivalent to having
$$\log\omega = u+\widetilde v,$$
where $u$ and $v$ are real-valued functions in $L^\infty(\T)$, with
$$\big\|v\big\|_{L^\infty(\T)}<\frac\pi2;$$
$\widetilde v$ is the harmonic conjugate of $v$. We note that this time,
the Carleson ``squares'' $Q$ are tacitly assumed to include the adjacent 
boundary arcs on $\T$. We interpret the Helson-Szeg\"o theorem as saying
that part of the BMO$(\T)$ norm of $\log\omega$ is small. A similar argument 
was used in \cite{Hed} to show that Brennan's conjecture is equivalent to an
area Muckenhoupt condition on $|\varphi'|^q$, for suitable exponents $q$.
The space that corresponds to the subspace BMOA$(\D)$ of BMO$(\T)$ 
(consisting of all functions whose Poisson extensions to the interior are 
holomorphic) in the case when arc length is replaced by area measure is 
the {\sl Bloch space} ${\mathcal B}(\D)$ (see, for instance, \cite{HKZ}) 
of all holomorphic functions $f$ in $\D$ with
$$\|f\|_{\mathcal B}=\sup\Big\{\big(1-|z|^2\big)\,|f'(z)|:\,\,z\in\D\Big\}
<+\infty;$$
the above expression is known as the {\sl Bloch norm}. This means that, 
{\sl ideologically, we should hope to find some estimates of the Bloch norm 
of $\log\varphi'$ which would be more or less equivalent to the 
$\diff A_\alpha$-area Muckenhoupt condition} (\ref{eq:intro-5}).
We are of course groping in the dark here, as there 
is no known theorem of Helson-Szeg\"o type that would apply in the (weighted)
area measure case. In any case, the conjectured property of 
${\mathrm B}_\schlicht(t)$ that this function is even near the origin, plus 
the related stronger rotational invariance suggested by Binder, lends 
credence to the idea that a study of (\ref{eq:intro-5}) is the same as 
studying the integral means spectral function ${\mathrm B}_\schlicht(\tau)$, 
at least for $\tau\in\C$ near the origin. If a function 
$f\in{\mathcal B}(\D)$ has sufficiently small Bloch norm, then it can be 
shown that $e^f$ belongs to any fixed Bergman space $\calH_\alpha(\D)$  
($-1<\alpha<+\infty$), with good control of the norm. In addition, it 
is also true  that $\log\varphi'\in{\mathcal B}(\D)$ for 
$\varphi\in\schlicht$; this is an 
easy consequence of (\ref{eq:Bieberbach}). The problem is that there is a 
genuine gap between the constants for the necessary and the sufficient 
conditions, and the only way to bridge that gap is to find an appropriate 
substitute for the Bloch norm as defined above. 
In \cite{Hed}, it was suggested by the first-named author, Hedenmalm, that 
spectral properties of a Volterra-type operator associated with 
$\log\varphi'$ should be relevant for the problem at hand; inspiration for 
this came from conversations with Alexandru Aleman. Essentially, this amounts 
to studying the multiplier properties of $\varphi''/\varphi'$. Then the 
second-named author, Shimorin, found that the multiplier norm of the 
Schwarzian derivative from the space $\calH_\alpha(\D)$ to 
$\calH_{\alpha+4}(\D)$ could be estimated effectively by using the area 
methods directly rather than going via the classical pointwise estimate
\begin{equation}
\left|\frac{\varphi'''(z)}{\varphi'(z)}-
\frac32\left[\frac{\varphi''(z)}{\varphi'(z)}\right]^2\right|
\le\frac6{(1-|z|^2)^2},\qquad z\in\D,
\label{eq:Bieberbach-2}
\end{equation}
and that this led to a better estimate of ${\rm B}_\schlicht(-1)$ and
${\rm B}_\schlicht(-2)$ than what was previously known. We should mention
that (\ref{eq:Bieberbach-2}) also expresses in a way that $\log\varphi'\in
{\mathcal B}(\D)$, and that the multiplier norm estimate implies
an estimate of the spectral radius of a Volterra-type operator associated
with the Schwarzian derivative. Shimorin's work suggests that the multiplier 
norm of the derivative of $\log\varphi'$ from $\calH_\alpha(\D)$ to 
$\calH_{\alpha+2}(\D)$ is a more appropriate way to measure the size of 
$\log\varphi'$ than applying the usual Bloch norm. Then, by dissecting 
a theorem by Prawitz, which generalizes the Gr\"onwall area theorem, we 
found a collection of estimates of multiplier norm type, parametrized by 
a real parameter $\theta$, $0<\theta\le1$. Generally speaking,
these estimates were the result of the application of the diagonal 
restriction operator on the bidisk $\D^2$ and the use of sharp constants 
in norm estimates.
By adding higher order terms corresponding to the multiplicity of the
zero along the diagonal, we found an estimate that was in fact an equality
for all full mappings $\varphi$. Unfortunately, the vast majority of these 
additional terms carry information of which it is, generally speaking, hard 
to make  
{\sl effective} use as regards the study of integral means spectra.
The details of the method are presented in Sections \ref{areathm}, 
\ref{se2}, and \ref{se4}.

\section{Area theorem type estimates}\label{areathm}

\noindent{\bf The theorem of Prawitz.} Our point of departure is a 
theorem of Prawitz, which generalizes Gr\"onwall's famous area theorem.

\begin{thm}
Let $\varphi\in\schlicht$. Then, for $0<\theta\le1$, we have
\begin{equation*}
%  \label{Prawitzineq}
  \int_\D\left|\varphi'(z)\,\left(\frac{z}{\varphi(z)}\right)^{\theta+1}
  -1\right|^2\frac{\diff A(z)}{|z|^{2\theta+2}}\le\frac1{\theta}\,,
\end{equation*}
with equality precisely for the full mappings $\varphi$. 
\label{thm-1}
\end{thm}

\begin{proof}
The inequality follows from a classical result of Prawitz, see
\cite[p. 13]{Milin} (the inequality in \cite{Milin}
is formulated for functions of the class $\Sigma$, but a standard passage
from $\Sigma$ to $\schlicht$ leads to the above inequality). The fact
that we have an equality precisely for the full mappings is a part of Prawitz'
theorem.
\end{proof}

In Theorem \ref{thm-1}, 
$$\left(\frac{z}{\varphi(z)}\right)^{\theta+1}=
\exp\left((\theta+1)\,\log\frac{z}{\varphi(z)}\right),$$
where the logarithm expression is determined uniquely by the requirements
that it be holomorphic in $\D$ and that it assume the value $0$ at $z=0$.
\medskip

\noindent{\bf A two-variable version of Prawitz' theorem.} 
We shall try to move the special point $z=0$ about in the disk, by the
following procedure. We start with a given $\varphi\in{\schlicht}$, and put
$$\psi(\zeta)=
\frac{\varphi\left(\frac{\zeta+w}{1+\bar w\zeta}\right)
-\varphi(w)}{(1-|w|^2)\,\varphi'(w)},\qquad\zeta\in\D,$$
for fixed $w\in\D$, which then is another element of $\schlicht$. 
Now, we
insert this $\psi$ in place of $\varphi$ in Theorem \ref{thm-1},
$$\int_\D\left|\frac1{\varphi'(w)}\,(1+\bar w \zeta)^{-2}\,
\varphi'\left(\frac{\zeta+w}{1+\bar w\zeta}\right)
\left(\frac{(1-|w|^2)\varphi'(w)\,\zeta}
{\varphi\left(\frac{\zeta+w}{1+\bar w\zeta}\right)
-\varphi(w)}\right)^{\theta+1}
-1\right|^2\frac{\diff A(\zeta)}{|\zeta|^{2\theta+2}}
\le\frac1{\theta}\,,$$ 
and we
make the change of variables
$$z=\frac{\zeta+w}{1+\bar w\zeta}\quad\Longleftrightarrow
\quad \zeta=\frac{z-w}{1-\bar w z}$$
in the integral. After simplification, we then obtain
\begin{equation}
\int_\D\left|\frac{\varphi'(z)}{\varphi'(w)}\,
\left(\frac{\varphi'(w)\,(z-w)}
{\varphi(z)-\varphi(w)}\right)^{\theta+1}
-\left(\frac{1-|w|^2}{1-\bar w z}\right)^{1-\theta}\right|^2
\frac{\diff A(z)}{|z-w|^{2\theta+2}}\le
\frac1{\theta}\,(1-|w|^2)^{-2\theta},
\label{eq-1}
\end{equation}
valid for all $\theta$ in the interval $0<\theta\le1$. 
This inequality is basic for our analysis. We shall write it in a slightly 
different form. Let
$$\Phi_\theta(z,w)=\frac{1}{z-w}\,\left\{
\frac{\varphi'(z)}{\varphi'(w)}\,
\left(\frac{\varphi'(w)\,(z-w)}
{\varphi(z)-\varphi(w)}\right)^{\theta+1}
-1\right\},\qquad (z,w)\in\D^2,\,\,\,z\neq w,$$
and
$$L_\theta(z,w)=\frac{1}{z-w}\,\left\{1-\left(
\frac{1-|w|^2}{1-\bar w z}\right)^{1-\theta}\right\},\qquad 
(z,w)\in\D^2,\,\,\,z\neq w.$$
We note that $\Phi_\theta$ extends analytically to the whole bidisk 
$\D^2$, and that its diagonal restriction is
$$\Phi_\theta(z,z)=
\frac{1-\theta}2\,\frac{\varphi''(z)}{\varphi'(z)}.$$
The function $L_\theta$ extends real analytically to $\D^2$.  
In view of (\ref{eq-1}), we have the following.

\begin{thm}
Fix $\theta$, $0<\theta\le1$, and let $\varphi\in\schlicht$ be arbitrary.
Then, for all $w\in\D$,
\begin{equation*}
\int_\D\Big|\Phi_\theta(z,w)+L_\theta(z,w)\Big|^2
\frac{\diff A(z)}{|z-w|^{2\theta}}
\le\frac1{\theta}\,(1-|w|^2)^{-2\theta},
\end{equation*}
with equality if and only if $\varphi$ is a full mapping.
\label{thm-2}
\end{thm}

\section{Bergman spaces in the bidisk}\label{se2}

For $-\infty<\alpha,\secpar<+\infty$, we consider the Hilbert space 
${\calL}_{\alpha,\secpar}(\D^2)$ of all Lebesgue measurable functions on 
the bidisk $\D^2$ (modulo null functions), subject to the norm boundedness 
condition
$$\|f\|_{\alpha,\secpar}=\left(\int_{\D}\int_\D |f(z,w)|^2|z-w|^{2\secpar}
\diff A(z)\,\diff A_\alpha(w)\right)^{1/2}<+\infty,$$
where $\diff A_\alpha$ is as in (\ref{eq-areael}).
 We also need the closed subspace 
${\calH}_{\alpha,\secpar}(\D^2)$ of ${\calL}_{\alpha,\secpar}(\D^2)$ that
consists of functions holomorphic in $\D^2$. {\sl The space 
${\calH}_{\alpha,\secpar}(\D^2)$ is trivial unless $-1<\alpha<+\infty$}. 
The reproducing kernel for the 
space ${\calH}_{\alpha,\secpar}(\D^2)$ will be denoted by
$$P_{\alpha,\secpar}\big((z,w);(z',w')\big),\qquad (z,w),\,(z',w')\in\D^2;$$
it is holomorphic in $(z,w)$, and anti-holomorphic in $(z',w')$. It is 
defined by the reproducing property
$$f(z,w)=\int_\D\int_\D P_{\alpha,\secpar}\big((z,w);(z',w')\big)\,
f(z',w')\,|z'-w'|^{2\secpar}\diff A(z')\,\diff A_\alpha(w'),$$
for all $(z,w)\in\D^2$ and $f\in{\calH}_{\alpha,\secpar}(\D^2)$. 
In case $\secpar=0$, it is given by the explicit formula
$$P_{\alpha,0}\big((z,w);(z',w')\big)=\frac{1}{(1-z\bar z')^2
(1-w\bar w')^{\alpha+2}},\qquad (z,w),\,(z',w')\in\D^2.$$
Associated with a kernel $T=T_{\alpha,\secpar}$ of the variables 
$((z,w);(z',w'))\in\D^2\times\D^2$, we have an operator on 
${\calL}_{\alpha,\secpar}(\D^2)$ defined by
$$\Top_{\alpha,\secpar}f(z,w)=\int_\D\int_\D 
T_{\alpha,\secpar}\big((z,w);(z',w')\big)\,f(z',w')\,|z'-w'|^{2\secpar}
\,\diff A(z')\,\diff A_\alpha(w'),$$
for $(z,w)\in\D^2$, which is going to be bounded in all cases we shall 
consider. For instance, associated with the kernel $P_{\alpha,\secpar}$ is 
the operator $\Pop_{\alpha,\secpar}$ which effects the orthogonal projection 
${\calL}_{\alpha,\secpar}(\D^2)\to{\calH}_{\alpha,\secpar}(\D^2)$.  

Let $N=0,1,2,3,\ldots$ be a nonnegative integer, and consider the closed 
subspace ${\calH}_{\alpha,\secpar;N}(\D^2)$ of 
${\calH}_{\alpha,\secpar}(\D^2)$ consisting of functions with
$$f(z,w)=O\big(|z-w|^N)$$
near the diagonal. These functions vanish up to degree $N$ along the diagonal,
and are holomorphically divisible by $(z-w)^N$. For $N=0$, we have
$${\calH}_{\alpha,\secpar;0}(\D^2)={\calH}_{\alpha,\secpar}(\D^2);$$
more generally, for $N=1,2,3,\ldots$,
$${\calH}_{\alpha,\secpar;N}(\D^2)={\calH}_{\alpha,\secpar}(\D^2)\qquad
\text{if}\quad -\infty<\beta+N\le0.$$
Being a closed subspace of the Hilbert space ${\calH}_{\alpha,\secpar}(\D^2)$,
the subspace  ${\calH}_{\alpha,\secpar;N}(\D^2)$ has a reproducing kernel 
function, denoted
$$P_{\alpha,\secpar;N}\big((z,w);(z',w')\big),\qquad (z,w),\,(z',w')\in\D^2.$$
Associated to the kernel is the orthogonal projection 
$$\Pop_{\alpha,\secpar;N}:\,{\calL}_{\alpha,\secpar}(\D^2)\to
{\calH}_{\alpha,\secpar;N}(\D^2).$$
The following is an important observation.

\begin{prop}
For $-1<\alpha,\secpar<+\infty$, we have
\begin{equation*}
P_{\alpha,\secpar;N}\big((z,w);(z',w')\big)=(z-w)^N(\bar z'-\bar w')^N
\, P_{\alpha,\secpar+N}\big((z,w);(z',w')\big),
\end{equation*}
for $(z,w),\,(z',w')\in\D^2$.
\label{prop-0.1}
\end{prop}

\begin{proof}
We note that multiplication by $(z-w)^N$ is an isometric isomorphism
$${\calH}_{\alpha,\secpar+N}(\D^2)\to{\calH}_{\alpha,\secpar;N}(\D^2);$$
from this, the conclusion is immediate.
\end{proof}

For $N=0,1,2,3,\ldots$, consider the Hilbert space 
$$\calI_{\alpha,\secpar;N}(\D^2)={\calH}_{\alpha,\secpar;N}(\D^2)\ominus
{\calH}_{\alpha,\secpar;N+1}(\D^2).$$
Its reproducing kernel has the form
\begin{equation}
Q_{\alpha,\secpar;N}\big((z,w);(z',w')\big)=
P_{\alpha,\secpar;N}\big((z,w);(z',w')\big)-
P_{\alpha,\secpar;N+1}\big((z,w);(z',w')\big),
\label{eq-5}
\end{equation}
and the associated operator projects orthogonally 
$$\Qop_{\alpha,\secpar;N}:\,{\calL}_{\alpha,\secpar}(\D^2)
\to{\calI}_{\alpha,\secpar;N}(\D^2).$$
We write $Q_{\alpha,\secpar}$ for the special kernel $Q_{\alpha,\secpar;0}$.
It then follows from Proposition \ref{prop-0.1} that
\begin{equation}
Q_{\alpha,\secpar;N}\big((z,w);(z',w')\big)=(z-w)^N(\bar z'-\bar w')^N
\, Q_{\alpha,\secpar+N}\big((z,w);(z',w')\big).
\label{eq-5.5}
\end{equation}

  The fact that the only function that vanishes to
an infinite degree along the diagonal is the zero function implies the
orthogonal decomposition 
$$ {\calH}_{\alpha,\secpar}(\D^2)=\bigoplus_{N=0}^{+\infty}
{\calI}_{\alpha,\secpar;N}(\D^2). $$
As a consequence, we have the decomposition of the kernel
\begin{multline}
\label{decompPalphabeta}
P_{\alpha,\secpar}\big((z,w);(z',w')\big)=
\sum_{N=0}^{+\infty}
Q_{\alpha,\secpar;N}\big((z,w);(z',w')\big)\\
=\sum_{N=0}^{+\infty}
(z-w)^N(\bar z'-\bar w')^N\, Q_{\alpha,\secpar+N}\big((z,w);(z',w')\big).
\end{multline}
and the norm decomposition
\begin{equation}
\label{normdecomp}
\big\|\Pop_{\alpha,\secpar}\,f\big\|^2_{\alpha,\secpar}=
\sum_{N=0}^{+\infty}\big\|\Qop_{\alpha,\secpar;N}\,f\big\|^2_{\alpha,\secpar},
\qquad f\in{\calL}_{\alpha,\secpar}(\D^2). 
\end{equation}

   There are some natural families  of unitary operators acting 
in spaces ${\calH}_{\alpha,\secpar}(\D^2)$. First, we can perform 
simultaneous rotations of variables $z$ and $w$: 
$$R_\theta[f](z,w)=f(e^{i\theta}z,e^{i\theta}w);
\quad f\in {\calH}_{\alpha,\secpar}(\D^2); \quad \theta\in\R. $$
The next family of unitary operators is given by the lemma below.

\begin{lemma}
For each $\lambda\in\D$, the operator
$$\Uop_\lambda [f](z,w)=\frac{(1-|\lambda|^2)^{\alpha/2+\secpar+2}}
{(1-\bar\lambda z)^{\secpar+2}(1-\bar\lambda w)^{\alpha+\secpar+2}}\,\,
f\left(\frac{\lambda-z}{1-\bar\lambda z},
\frac{\lambda-w}{1-\bar\lambda w}\right)$$
is unitary on $\calH_{\alpha,\secpar}(\D^2)$, and its square is the identity:
$\Uop_\lambda^2 [f]=f$ for all $f\in\calH_{\alpha,\secpar}(\D^2)$.
\label{lm-5}
\end{lemma}

\begin{proof} This amounts to an elementary change of variables calculation.
\end{proof}

In fact, if both $\alpha$ and $\secpar$ are even integers, then 
for each M\"obius automorphism $\psi$ of the disk $\D$, 
one can  define the operator $\Uop_\psi$: 
$$\Uop_\psi[f](z,w)=f\big(\psi(z),\psi(w)\big)\cdot
\big(\psi'(z) \big)^{1+\secpar/2}\cdot
\big(\psi'(w)\big)^{1+\alpha/2+\beta/2}. $$
Then all operators $\Uop_\psi$ are unitary in $\calH_{\alpha,\secpar}(\D^2)$
and the map $\psi\mapsto \Uop_\psi$ is a unitary representation of the 
group of M\"obius automorphisms of $\D$.

 We proceed  by analyzing the reproducing kernel 
$P_{\alpha,\secpar}$
along the diagonal.

\begin{lemma} Fix $-1<\alpha,\secpar<+\infty$.
We then have 
$$P_{\alpha,\secpar}\big((z,w);(z',z')\big)=
Q_{\alpha,\secpar}\big((z,w);(z',z')\big)=\frac{\thpar(\alpha,\secpar)}
{(1-z\bar z')^{\secpar+2}(1-w\bar z')^{\alpha+\secpar+2}},$$
where the constant $\thpar(\alpha,\secpar)$ is given by
$$\frac1{\thpar(\alpha,\secpar)}=\int_\D\int_\D|z-w|^{2\secpar}
\diff A(z)\,\diff A_\alpha(w).$$
\label{lm-1}
\end{lemma}

\begin{proof}
We note first that the fact that rotation  operators $R_\theta$ are
unitary in $\calH_{\alpha,\secpar}(\D^2)$ implies that 
$$P_{\alpha,\secpar}\big((e^{i\theta}z,e^{i\theta }w );(0,0)\big)= 
P_{\alpha,\secpar}\big((z,w );(0,0)\big). $$
Now, we observe that the only functions analytic in $\D^2$ and having this 
property are the constant functions, which follows at once by considering 
double power series expansions. Hence, 
$P_{\alpha,\secpar}\big((z,w);(0,0)\big)$ is constant in $(z,w)$, and we
write
\begin{equation}
\thpar(\alpha,\secpar)=P_{\alpha,\secpar}\big((z,w);(0,0)\big)
\label{eq-repr-0}
\end{equation}
for this constant. The above integral formula for $\thpar(\alpha,\secpar)$ 
follows from the reproducing property of the kernel
$P_{\alpha,\secpar}\big((\cdot,\cdot );(0,0)\big)$ applied to the constant 
function $1$. 

Now, let $\lambda\in \D$. We pick $f\in\calH_{\alpha,\secpar}(\D^2)$, and 
note that in view of (\ref{eq-repr-0}) and Lemma \ref{lm-5},
\begin{multline*}
\big(1-|\lambda|^2\big)^{\alpha/2+\beta+2}\,f(\lambda,\lambda)
=\Uop_\lambda [f](0,0)=\thpar(\alpha,\secpar)\,
\big\langle\Uop_\lambda^2 [f],\Uop_\lambda[1]\big\rangle_{\alpha,\secpar}
\\
=\thpar(\alpha,\secpar)\,
\big\langle f,\Uop_\lambda[1]\big\rangle_{\alpha,\secpar}.
\end{multline*}
This formula expresses the reproducing identity
at the diagonal point $(\lambda,\lambda)$, which shows that 
$$P_{\alpha,\secpar}\big((z,w);(\lambda,\lambda)\big)=
\thpar(\alpha, \secpar)(1-|\lambda|^2)^{-\alpha/2-\beta-2}
\Uop_\lambda[1](z,w);$$
after some simplification, this gives the desired expression. 
\end{proof}

In view of Lemma \ref{lm-1}, 
$$P_{\alpha,\secpar}\big((z,z);(z',z')\big)=\frac{\thpar(\alpha,\secpar)}
{(1-z\bar z')^{\alpha+2\secpar+4}},$$
which we identify as the reproducing kernel for the 
Hilbert space coinciding as a set with the space 
$\calH_{\alpha+2\secpar+2}(\D)$ 
from the introduction and supplied with the norm 
$$\|f\|^2
=\frac1{\thpar(\alpha,\secpar)}\int_\D|f(z)|^2\diff A_{\alpha+2\secpar+2}(z)
=\frac 1{\thpar(\alpha, \secpar)}\|f\|^2_{\alpha+2\secpar+2}.$$
Let $\oslash$ denote the operation of taking the diagonal restriction:
$$(\oslash f)(z)=f(z,z),\qquad z\in\D.$$
In view of the general theory of reproducing kernels (see \cite{Aronsz}
and \cite{Saitoh}),  we have the sharp estimate
\begin{equation}
\frac1{\thpar(\alpha,\secpar)}\,\|\!\oslash f\|_{\alpha+2\secpar+2}^2\le
\|f\|^2_{\alpha,\secpar},\qquad f\in \calH_{\alpha,\secpar}(\D^2).
\label{eq-diagest}
\end{equation}
In fact, we can even determine the corresponding norm identity.

\begin{lemma}
We have the equality of norms
$$\frac1{\thpar(\alpha,\secpar)}\,\big\|\!\oslash 
f\big\|_{\alpha+2\secpar+2}^2=
\big\|\Qop_{\alpha,\secpar}f\big\|^2_{\alpha,\secpar},
\qquad f\in \calH_{\alpha,\secpar}(\D^2).$$
\label{lm-2}
\end{lemma}

\begin{proof}
The analysis of reproducing kernel functions that leads up to the estimate
(\ref{eq-diagest}) also shows that to each 
$f\in \calH_{\alpha,\secpar}(\D^2)$ there exists a 
$g\in \calH_{\alpha,\secpar}(\D^2)$ such that $\oslash g=\oslash f$ and
$$\frac1{\thpar(\alpha,\secpar)}\,\|\!\oslash f\|_{\alpha+2\secpar+2}^2=
\|g\|^2_{\alpha,\secpar},\qquad f\in \calH_{\alpha,\secpar}(\D^2).$$
We decompose this $g$ as follows:
$$g=\Qop_{\alpha,\secpar}f+\big(g-\Qop_{\alpha,\secpar}f\big)\in
\calI_{\alpha,\secpar;0}(\D^2)+\calH_{\alpha,\secpar;1}(\D^2).$$
As this decomposition is orthogonal, we get
$$\big\|\Qop_{\alpha,\secpar}f\big\|_{\alpha,\secpar}^2\le
\big\|\Qop_{\alpha,\secpar}f\big\|_{\alpha,\secpar}^2+
\big\|g-\Qop_{\alpha,\secpar}f\big\|_{\alpha,\secpar}^2=
\Big\|\Qop_{\alpha,\secpar}f+
\big(g-\Qop_{\alpha,\secpar}f\big)\Big\|_{\alpha,\secpar}^2=
\|g\|^2_{\alpha,\secpar}.$$
The assertion now follows from the above estimates together with
(\ref{eq-diagest}).
\end{proof}

The constant $\thpar(\alpha,\secpar)$ can be evaluated explicitly. 

\begin{lemma} Fix $-1<\alpha,\secpar<+\infty$.
Then 
\begin{equation*}
%\label{aalphabeta}
\frac1{\thpar(\alpha,\secpar)}=\int_\D\int_\D|z-w|^{2\secpar}\,\diff A(z)
\,\diff A_\alpha(w)
=\frac1{1+\secpar}\,\frac{\Gamma(\alpha+2)\,\Gamma(\alpha+2\secpar+3)}
{\Gamma(\alpha+\secpar+2)\,\Gamma(\alpha+\secpar+3)}.
\end{equation*}
\label{lm-4}
\end{lemma}

\begin{proof}
We perform the change of variables
$$\zeta=\frac{w-z}{1-\bar w\,z},\qquad z=\frac{w-\zeta}{1-\bar w\,\zeta},$$
and replace the pair $(z,w)$ by $(\zeta,w)$. The result is, after
simplification,
\begin{multline*}
\frac1{\thpar(\alpha,\secpar)}=\frac{\alpha+1}{(1+\secpar)(\alpha+2\secpar+3)}
\sum_{n=0}^{+\infty}\frac{(\secpar+2)_n(\secpar+1)_n}
{n!\,(\alpha+2\secpar+4)_n}\\
=\frac{\alpha+1}{(1+\secpar)(\alpha+2\secpar+3)}\,
{}_2F_1\big(\secpar+2,\secpar+1;\alpha+2\secpar+4;1\big),
\end{multline*}
where ${}_2F_1$ denotes Gauss' hypergeometric function. Here, we use the 
standard Pochhammer notation
$$(a)_n=a(a+1)(a+2)\ldots(a+n-1).$$
The assertion now follows from the well-known identity
\begin{equation}
\label{Gausshypergeom}
{}_2F_1(a,b;c;1)=\frac{\Gamma(c)\,\Gamma(c-a-b)}{\Gamma(c-a)\Gamma(c-b)}.
\end{equation}
The proof is complete.
\end{proof}

\begin{remark} \rm
It follows from Lemma \ref{lm-4} that
$$\frac{\thpar(\alpha,\secpar+n)}{\thpar(\alpha,\secpar)}=
\frac{n+1+\secpar}{1+\secpar}
\,\frac{(\alpha+\secpar+2)_n(\alpha+\secpar+3)_n}{(\alpha+2\secpar+3)_{2n}},
\qquad n=1,2,3,\ldots.$$
\label{rm-1}
\end{remark}

We obtain an integral representation of the kernel $Q_{\alpha,\secpar}$.

\begin{lemma} Fix $-1<\alpha,\secpar<+\infty$. 
The kernel $Q_{\alpha,\secpar}$ is given by the integral formula
$$Q_{\alpha,\secpar}\big((z,w);(z',w')\big)=\thpar(\alpha,\secpar)\int_\D
\frac{\diff A_{\alpha+2\secpar+2}(\xi)}{(1-\bar\xi z)^{\secpar+2}
(1-\bar\xi w)^{\alpha+\secpar+2}(1-\xi\bar z')^{\secpar+2}
(1-\xi\bar w')^{\alpha+\secpar+2}},$$
for $(z,w),(z',w')\in\D^2$.
\label{lm-3}
\end{lemma}

\begin{proof}
It is enough to establish that if $\widetilde Q_{\alpha,\secpar}$ denotes
the kernel defined by the above integral formula, then it coincides with
the reproducing kernel function for the space $\calI_{\alpha,\secpar;0}
(\D^2)=\calH_{\alpha,\secpar}(\D^2)\ominus\calH_{\alpha,\secpar;1}(\D^2)$.
To this end, we first check that for each individual point $(z',w')\in\D^2$,
the function
\begin{multline*}
(z,w)\mapsto\widetilde Q_{\alpha,\secpar}\big((z,w);(z',w')\big)\\
=\thpar(\alpha,\secpar)\int_\D
\frac{\diff A_{\alpha+2\secpar+2}(\xi)}{(1-\bar\xi z)^{\secpar+2}
(1-\bar\xi w)^{\alpha+\secpar+2}(1-\xi\bar z')^{\secpar+2}
(1-\xi\bar w')^{\alpha+\secpar+2}},
\end{multline*}
belongs to 
$\calH_{\alpha,\secpar}(\D^2)\ominus\calH_{\alpha,\secpar;1}(\D^2)$. As
a first step, we see that if we use the methods of Chapter 1 in \cite{HKZ},
we can show that this function belongs to $\calL_{\alpha,\secpar}(\D^2)$,
and then, by inspection, it is also analytic in $\D^2$, and hence 
an element of $\calH_{\alpha,\secpar}(\D^2)$. To prove that it is 
orthogonal to $\calH_{\alpha,\secpar;1}(\D^2)$, we note that each ``term''
\begin{equation*}
(z,w)\mapsto\frac{1}{(1-\bar\xi z)^{\secpar+2}
(1-\bar\xi w)^{\alpha+\secpar+2}(1-\xi\bar z')^{\secpar+2}
(1-\xi\bar w')^{\alpha+\secpar+2}},
\end{equation*}
is a multiple of the element that achieves the point evaluation at the
diagonal point $(\xi,\xi)$, and therefore it is orthogonal to 
the subspace $\calH_{\alpha,\secpar;1}(\D^2)$, as these functions vanish 
at all diagonal points.

   Now, 
we  see, by inspection, that
\begin{equation*}
\widetilde Q_{\alpha,\secpar}\big((z,z);
(z',w')\big)=\frac{\thpar(\alpha,\secpar)}
{(1-z\bar z')^{\secpar+2}(1-z\bar w')^{\alpha+\secpar+2}};
\end{equation*}
this follows from 
 the reproducing property of
the well-known kernel function in the space $\calH_{\alpha+2\secpar+2}(\D)$.
 We note that this is the same as $Q_{\alpha,\secpar}\big((z,z);
(z',w')\big)$, according to Lemma \ref{lm-1}.
And since functions from $\calI_{\alpha,\secpar;0}$ are 
uniquely determined by their diagonal restrictions, we obtain 
$\widetilde
Q_{\alpha,\secpar}=Q_{\alpha,\secpar}$.
\end{proof}

\begin{prop}
Fix $-1<\alpha,\secpar<+\infty$. Then, for $N=0,1,2,3,\ldots$, we have
\begin{equation*}
\big\|\Qop_{\alpha,\secpar;N}\,f\big\|^2_{\alpha,\secpar}=
\frac1{\thpar(\alpha,\secpar+N)}\,
\left\|\,\oslash \left[\frac{\Pop_{\alpha,\secpar;N}\,f(z,w)}{(z-w)^N}
\right]\right\|^2_{\alpha+2\secpar+2N+2},\qquad 
f\in{\calL}_{\alpha,\secpar}(\D^2).
\end{equation*}
\label{prop-7.11}
\end{prop}

\begin{proof}
This follows from a combination of Proposition \ref{prop-0.1} 
and Lemma \ref{lm-2}.
\end{proof}

In view of Lemma \ref{lm-1}, we have, for $z\in\D$,
\begin{multline}
\oslash \left[\frac{\Pop_{\alpha,\secpar;N}\,f(z,w)}{(z-w)^N}\right](z)=
\thpar(\alpha,\secpar+N)\\
\times\int_\D\int_\D
\frac{(\bar z'-\bar w')^N}
{(1-z\bar z')^{\secpar+N+2}(1-z\bar w')^{\alpha+\secpar+N+2}}
\,f(z',w')\,|z'-w'|^{2\secpar}\diff A(z')\,\diff A_\alpha(w').
\label{eq-8}
\end{multline}
 
We want to express this in terms of derivatives of order $N$ of $f$. 
To this end, we note that the series expansion in 
(\ref{decompPalphabeta}) leads to
\begin{equation}
\big[\bpartial_z^k\, P_{\alpha,\secpar}\big]\big((z,z);(z',w')\big)
=\sum_{n=0}^{k}n!\,\big(\bar z'-\bar w'\big)^n
\left(\begin{array}{c}k\\n\end{array}\right)
\,\Big[\bpartial_z^{k-n}\, 
Q_{\alpha,\secpar+n}\Big]\big((z,z);(z',w')\big),
\label{eq-9}
\end{equation}
where $\bpartial_z$ stands for the (partial) derivative with respect to $z$. 
Moreover, in view of Lemma \ref{lm-3},
\begin{multline*}
\bpartial_z^{k-n}\,Q_{\alpha,\secpar}\big((z,w);(z',w')\big)=
\thpar(\alpha,\secpar)\,(\secpar+2)_{k-n}\\
\times\int_\D
\frac{\bar\xi^{k-n}\,\diff A_{\alpha+2\secpar+2}(\xi)}
{(1-\bar\xi z)^{\secpar+k-n+2}
(1-\bar\xi w)^{\alpha+\secpar+2}(1-\xi\bar z')^{\secpar+2}
(1-\xi\bar w')^{\alpha+\secpar+2}},
\end{multline*}
which, when restricted to the diagonal, becomes
\begin{multline*}
\big[\bpartial_z^{k-n}\,Q_{\alpha,\secpar}\big]
\big((z,z);(z',w')\big)=
\thpar(\alpha,\secpar)\,(\secpar+2)_{k-n}\\
\times\int_\D
\frac{\bar\xi^{k-n}\,\diff A_{\alpha+2\secpar+2}(\xi)}
{(1-\bar\xi z)^{\alpha+2\secpar+k-n+4}(1-\xi\bar z')^{\secpar+2}
(1-\xi\bar w')^{\alpha+\secpar+2}}\\
=\frac{(\secpar+2)_{k-n}}{(\alpha+2
\secpar+4)_{k-n}}\,\,\,\bpartial_z^{k-n}\,\frac{\thpar(\alpha,\secpar)}
{(1-z\bar z')^{\secpar+2}(1-z\bar w')^{\alpha+\secpar+2}}.
%\label{eq-10}
\end{multline*}
By changing $\secpar$ to $\secpar+n$, we obtain, in view of Lemma \ref{lm-1},
that 
\begin{multline}
\big[\bpartial_z^{k-n}\,Q_{\alpha,\secpar+n}\big]
\big((z,z);(z',w')\big)\\=
\frac{(\secpar+n+2)_{k-n}}{(\alpha+2
\secpar+2n+4)_{k-n}}\,\,\,\bpartial_z^{k-n}\,\Big[P_{\alpha,\secpar+n}
\big((z,z);(z',w')\big)\Big].
\label{eq-10'}
\end{multline}
Now, applying (\ref{eq-9}) to a function $f\in\calL_{\alpha,\secpar}(\D^2)$, 
while taking (\ref{eq-10'}) into account, we find that
\begin{multline*}
\oslash\big[\bpartial_z^k\Pop_{\alpha,\secpar}f\big](z)=\sum_{n=0}^k
n!\,\left(\begin{array}{c}k\\n\end{array}\right)\,
\,\frac{(\secpar+n+2)_{k-n}}{(\alpha+2\secpar+2n+4)_{k-n}}
\,\,\bpartial_z^{k-n}\,
\oslash \left[\frac{\Pop_{\alpha,\secpar;n}\,f(z,w)}{(z-w)^n}\right](z).
\end{multline*}
We differentiate the above relation $N-k$ times with respect to $z$,
and obtain
\begin{multline}
\bpartial_z^{N-k}\oslash\big[\bpartial_z^k\Pop_{\alpha,\secpar}f\big](z)\\
=\sum_{n=0}^k
n!\,\left(\begin{array}{c}k\\n\end{array}\right)\,
\,\frac{(\secpar+n+2)_{k-n}}{(\alpha+2
\secpar+2n+4)_{k-n}}\,\,\,\bpartial_z^{N-n}\,\,
\oslash \left[\frac{\Pop_{\alpha,\secpar;n}\,f(z,w)}{(z-w)^n}\right](z).
\label{eq-10.2}
\end{multline}

We now formulate the desired relation.

\begin{prop} Fix $-1<\alpha,\beta<+\infty$. 
For each $N=0,1,2,3,\ldots$, we have 
\begin{equation*}
\oslash \left[\frac{\Pop_{\alpha,\secpar;N}\,f(z,w)}{(z-w)^N}\right]=
\sum_{k=0}^N a_{k,N}\,\,\bpartial_z^{N-k}\oslash\big[\bpartial_z^k 
\Pop_{\alpha,\secpar}f\big],
\end{equation*}
where 
$$a_{k,N}=\frac{(-1)^{N-k}}{k!(N-k)!}\,\,\,\frac{(\secpar+k+2)_{N-k}}
{(\alpha+2\secpar+N+k+3)_{N-k}}.$$
\label{prop-1.5}
\end{prop}

\begin{proof}
In view of (\ref{eq-10.2}), we should verify that
\begin{multline}
\sum_{k=0}^N a_{k,N}\,\,
\bpartial_z^{N-k}\oslash\big[\bpartial_z^k\Pop_{\alpha,\secpar}f\big](z)\\
=\sum_{k=0}^N\sum_{n=0}^k a_{k,N}\,\,
n!\,\left(\begin{array}{c}k\\n\end{array}\right)\,
\,\frac{(\secpar+n+2)_{k-n}}{(\alpha+2
\secpar+2n+4)_{k-n}}\,\,\bpartial_z^{N-n}\,
\oslash \left[\frac{\Pop_{\alpha,\secpar;n}\,f(z,w)}{(z-w)^n}\right](z)\\
=\oslash \left[\frac{\Pop_{\alpha,\secpar;N}\,f(z,w)}{(z-w)^N}\right](z),
\label{eq-10.3}
\end{multline}
where $a_{k,N}$ is as above. We realize that it is enough to show that
\begin{equation*}
\sum_{k=n}^N a_{k,N}\,\,
n!\,\left(\begin{array}{c}k\\n\end{array}\right)\,
\,\frac{(\secpar+n+2)_{k-n}}{(\alpha+2
\secpar+2n+4)_{k-n}}=\delta_{n,N},\qquad n=0,1,2,3,\ldots,N,
\end{equation*}
where the delta is the usual Kronecker symbol; as we implement the given 
values of the constants $a_{k,N}$, this amounts to
\begin{equation*}
\sum_{k=n}^N \frac{(-1)^{N-k}}{(k-n)!\,(N-k)!}\,
\frac{(\secpar+n+2)_{k-n}(\secpar+k+2)_{N-k}}{(\alpha+2
\secpar+2n+4)_{k-n}(\alpha+2\secpar+N+k+3)_{N-k}}=\delta_{n,N},
\end{equation*}
for $n=0,1,2,3,\ldots,N$. We quickly verify that this is correct for $n=N$.
To deal with smaller values of $n$, we first note that 
$$(\secpar+n+2)_{k-n}(\secpar+k+2)_{N-k}=(\secpar+n+2)_{N-n},$$
which is independent of $k$, so that we may factor it out, and reduce the 
problem to showing that
\begin{equation*}
\sum_{k=n}^N \frac{(-1)^{N-k}}{(k-n)!\,(N-k)!}\,
\frac{1}{(\alpha+2\secpar+2n+4)_{k-n}(\alpha+2\secpar+N+k+3)_{N-k}}
=0,
\end{equation*}
for $n=0,1,2,\ldots,N-1$. We compute that 
$$(\alpha+2\secpar+2n+4)_{k-n}(\alpha+2\secpar+N+k+3)_{N-k}=
\frac{(\alpha+2\secpar+2n+4)_{2N-2n-1}}{(\alpha+2\secpar+n+k+4)_{N-n-1}},$$
which reduces our task further to showing that
\begin{equation*}
\sum_{k=n}^N \frac{(-1)^{N-k}}{(k-n)!\,(N-k)!}\,
(\alpha+2\secpar+n+k+4)_{N-n-1}=0,
\end{equation*}
for $n=0,1,2,\ldots,N-1$.
We introduce the variables $j=k-n$ and $N'=N-n$, and rewrite the above:
\begin{equation*}
\sum_{j=0}^{N'} \frac{(-1)^{N'-j}}{j!\,(N'-j)!}\,
(\alpha+2\secpar+2n+j+4)_{N'-1}=0,
\end{equation*}
for $n=0,1,2,\ldots$ and $N'=1,2,3,\ldots$. Next, we consider the variable
$$\lambda=\alpha+2\secpar+2n+4,$$
which we shall think of as an independent variable, and we once more 
rewrite the above assertion:
\begin{equation*}
\sum_{j=0}^{N'}(-1)^{j}\,
\left(\begin{array}{c}N'\\j\end{array}\right)\,
(\lambda+j)_{N'-1}=0,
\end{equation*}
for $N'=1,2,3,\ldots$. The expression $q(\lambda)=(\lambda)_{N'-1}$ is a 
polynomial of degree $N'-1$ in $\lambda$, and
\begin{equation*}
\sum_{j=0}^{N'}(-1)^{j}\,
\left(\begin{array}{c}N'\\j\end{array}\right)\,
q(\lambda+j)
\end{equation*}
is an $N'$-th order iterated difference, which automatically produces $0$
on polynomials of degree less than $N'$. The assertion follows.
\end{proof}

We finally obtain an expansion of the norm in $\calH_{\alpha,\secpar}(\D^2)$
on the bidisk in terms of ``one-dimensional'' norms, taken over the unit disk,
analogous to the Taylor expansion along the diagonal.

\begin{cor}
For $f\in\calH_{\alpha,\secpar}(\D^2)$, we have the norm expansion
$$\|f\|_{\alpha,\beta}^2=
\sum_{N=0}^{+\infty}\frac1{\thpar(\alpha,\secpar+N)}\,
\,\bigg\|\sum_{k=0}^N
a_{k,N}\,\,\bpartial_z^{N-k}\oslash\big[\bpartial_z^k f\big]\bigg\|^2
_{\alpha+2\secpar+2N+2},$$
where the constants are as in Lemma \ref{lm-4} and Proposition \ref{prop-1.5}.
\label{cor-1.10}
\end{cor}

\begin{proof}
This results from a combination of (\ref{normdecomp})
and Propositions  \ref{prop-7.11} 
and \ref{prop-1.5}.
\end{proof}

\section{The main inequality}\label{se4}

\noindent{\bf Integration with respect to the second variable.}
 Fix $\theta$, $0<\theta\le1$, and let $\varphi\in{\schlicht}$ be arbitrary.
At times, the calculations below will be valid only for $0<\theta<1$, but
the validity for $\theta=1$ can usually be established easily by a simple
limit argument. By Theorem \ref{thm-2}, we have
\begin{equation}
\int_\D\Big|
\Phi_\theta(z,w)+L_\theta(z,w)\Big|^2
\frac{\diff A(z)}{|z-w|^{2\theta}}\le
\frac1{\theta}\,(1-|w|^2)^{-2\theta},
\label{eq-16}
\end{equation}
Let $g$ be a function that is holomorphic in $\D$. Then, in view
of (\ref{eq-16}), 
\begin{multline}
\int_\D\int_\D\Big|
\Phi_\theta(z,w)\,g(w)+L_\theta(z,w)\,g(w)\Big|^2\,\,|z-w|^{-2\theta}
\,\diff A(z)\,\diff A_\alpha(w)\\
\le
\frac1{\theta}\,\int_\D |g(w)|^2\,(1-|w|^2)^{-2\theta}\,\diff A_\alpha(w)=
\frac{\alpha+1}{\theta(\alpha-2\theta+1)}\|g\|^2_{\alpha-2\theta}
\label{eq-17}
\end{multline} 
(the last equality holds 
provided that 
$-1+2\theta<\alpha<+\infty$). 

In what follows, we assume that $g\in\calH_{\alpha-2\theta}(\D)$
and $-1<\alpha-2\theta<+\infty$. The left hand side of (\ref{eq-17})
expresses the square of the norm of the function 
$\Phi_\theta(z,w)\,g(w)+L_\theta(z,w)\,g(w)$
in the space $\calL_{\alpha,-\theta}(\D^2)$. It will be shown later that 
both terms of this sum belong to  
$\calL_{\alpha,-\theta}(\D^2)$ and hence one has the following decomposition:
\begin{multline}
\Phi_\theta(z,w)\,g(w)+L_\theta(z,w)\,g(w)\\
=\Big\{\Phi_\theta(z,w)\,g(w)+
\Pop_{\alpha,-\theta}\big[L_\theta(z,w)\,g(w)\big]\Big\}
+\Pop^\perp_{\alpha,-\theta}\big[L_\theta(z,w)\,g(w)\big],
\label{eq-10.55}
\end{multline}
with the corresponding decomposition of the norm
\begin{multline}
\big\|\Phi_\theta(z,w)\,g(w)+L_\theta(z,w)\,g(w)\big\|^2_{\alpha,-\theta}\\
=\big\|\Phi_\theta(z,w)\,g(w)+
\Pop_{\alpha,-\theta}\big[L_\theta(z,w)\,g(w)\big]\big\|^2_{\alpha,-\theta}
+\big\|\Pop^\perp_{\alpha,-\theta}
\big[L_\theta(z,w)\,g(w)\big]\big\|^2_{\alpha,-\theta}.
\label{main-1}
\end{multline}
Here, $\Pop^\perp_{\alpha,-\theta}$ is the projection complementary to 
$\Pop_{\alpha,-\theta}$: 
$$\Pop^\perp_{\alpha,-\theta}={\mathbf I}-\Pop_{\alpha,-\theta}
\ \mbox{in }\ \calL_{\alpha,-\theta}(\D^2),$$
where ${\mathbf I}$ stands for the identity operator. It follows that the 
inequality (\ref{eq-17}) assumes the form
\begin{multline}
\label{main-2}
\big\|\Phi_\theta(z,w)\,g(w)+
\Pop_{\alpha,-\theta}\big[L_\theta(z,w)\,g(w)\big]\big\|^2_{\alpha,-\theta}
  \\  
\le \frac{\alpha+1}{\theta(\alpha-2\theta+1)}\|g\|^2_{\alpha-2\theta}
-\big\|\Pop^\perp_{\alpha,-\theta}
\big[L_\theta(z,w)\,g(w)\big]\big\|^2_{\alpha,-\theta}.
\end{multline}
\medskip

\noindent{\bf The norm of a projected term.} We shall find an explicit
expression for the squared norm 
$$\big\|\Pop^\perp_{\alpha,-\theta}
\big[L_\theta(z,w)\,g(w)\big]\big\|^2_{\alpha,-\theta}.$$
 We first note that 
\begin{equation}
\label{main-3}
\big\|\Pop^\perp_{\alpha,-\theta}
\big[L_\theta(z,w)\,g(w)\big]\big\|^2_{\alpha,-\theta}= 
 \big\| L_\theta(z,w)\,g(w)\big\|^2_{\alpha,-\theta}
-\big\|\Pop_{\alpha,-\theta}
\big[L_\theta(z,w)\,g(w)\big]\big\|^2_{\alpha,-\theta}
\end{equation}
 We recall the classical definition of the Gauss hypergeometric function:
$${}_2F_1(a,b;c;x)=1+\sum_{n=1}^{+\infty}\frac{(a)_n(b)_n}{(c)_n\,n!}\,x^n,$$
where the series converges at least for complex $x\in\D$, unless we 
accidentally divide by zero in any of the terms.

\begin{lemma} 
For fixed $w\in\D$, we have the identity
\begin{equation*}
\int_\D\big|L_\theta(z,w)\big|^2\,\,|z-w|^{-2\theta}
\,\diff A(z)=\frac1{\theta}\,\Big[1-{}_2F_1
\big(1-\theta,-\theta;1;|w|^2\big)\Big]\,
\big(1-|w|^2\big)^{-2\theta}.
\end{equation*}
\label{lm-5'}
\end{lemma}

\begin{proof}
We make the change of variables
$$z=\frac{w-\zeta}{1-\bar w\zeta}, \qquad \zeta=\frac{w-z}{1-\bar w z},$$
which gives
\begin{equation*}
\int_\D\big|L_\theta(z,w)\big|^2\,\,|z-w|^{-2\theta}
\,\diff A(z)=\big(1-|w|^2\big)^{-2\theta}
\int_\D\left|\frac{1-(1-\bar w\zeta)^{\theta-1}}{\zeta}\right|^2
\,|\zeta|^{-2\theta}\,\diff A(\zeta).
\end{equation*}
We expand the power appearing in the integrand on the right hand side as
a Taylor series, and use that $z^j$ and $z^k$ are orthogonal in a radially
weighted Bergman space whenever $j\neq k$. The expression involving the
Gauss hypergeometric function then results from this.  
\end{proof}

\begin{lemma} 
For $w\in\D$, we have 
$${}_2F_1\big(1-\theta,-\theta;1;|w|^2\big)\ge{}_2F_1
\big(1-\theta,-\theta;1;1\big)
=\frac{\Gamma(2\theta+1)}{2[\Gamma(\theta+1)]^2}.$$
\label{lm-6}
\end{lemma}

\begin{proof}
The inequality follows if we see that the coefficients of the Taylor
series for ${}_2F_1(1-\theta,-\theta;1;x)$ are all negative except for the 
first one. The evaluation of 
\linebreak 
${}_2F_1(1-\theta,-\theta;1;1)$ is classical 
(see any book on special functions).
\end{proof}

Combining these two lemmas, we obtain the following.

\begin{prop}
For $g\in\calH_{\alpha-2\theta}(\D)$, we have
\begin{multline*}
\int_\D\big|L_\theta(z,w)\,g(w)\big|^2\,\,|z-w|^{-2\theta}
\,\diff A(z)\,\diff A_\alpha(w)\\
=\frac{\alpha+1}{\theta\,(\alpha-2\theta+1)}\int_\D\Big[1-{}_2F_1
\big(1-\theta,-\theta;1;|w|^2\big)\Big]\,|g(w)|^2\,
\diff A_{\alpha-2\theta}(w)\\
\le\frac{\alpha+1}{\theta\,(\alpha-2\theta+1)}\,
\left[1-\frac{\Gamma(2\theta+1)}{2[\Gamma(\theta+1)]^2}\right]\,
\|g\|^2_{\alpha-2\theta}.
\end{multline*}
\label{prop-1.23} 
\end{prop}
In particular, we see that the function $L_\theta(z,w)g(w)$ 
is in the space $\calL_{\alpha,-\theta}(\D^2)$. For later use, we need the 
following representation of the square of its norm: 
\begin{multline}
\label{main-4}
\big\|L_\theta(z,w)g(w)\big\|^2_{\alpha,-\theta}= 
\frac{\alpha+1}{\theta(\alpha-2\theta+1)}
\left(1-\frac{\Gamma(2\theta+1)}{2[\Gamma(\theta+1)]^2}\right)
\|g\|^2_{\alpha-2\theta}+  
\frac{\alpha+1}{\theta(\alpha-2\theta+1)}\\
\times
\int_\D \Big[{}_2F_1(1-\theta,-\theta;1;1)-
{}_2F_1(1-\theta,-\theta;1;|w|^2)\Big]\,|g(w)|^2\,
\diff A_{\alpha-2\theta}(w). 
\end{multline}

Only the first term of this sum is essential for our purposes, and 
the second may be considered as a contribution of ``higher order terms''. 
This is made explicit in the following lemma.

\begin{lemma}
\label{main-estimate1}
There exists a positive constant $C_1=C_1(\alpha, \theta)$ depending only on 
$\alpha$ and $\theta$ such that 
\begin{equation*}
0\le\int_\D \Big[{}_2F_1(1-\theta,-\theta;1;|w|^2)-
{}_2F_1(1-\theta,-\theta;1;1)\Big]\,|g(w)|^2\,\diff A_{\alpha-2\theta}(w)
\le C_1\, \|g\|^2_{\alpha-\theta}. 
\end{equation*}
\label{lm-7}
\end{lemma}

\begin{proof}
We use the inequality 
$$1-x^n\le n^\theta (1-x)^\theta, \qquad 0\le x\le1,$$
and the well-known asymptotics of the Pochhammer symbol
$$\frac{(1-\theta)_n}{n!}\sim \frac{n^{-\theta}}{\Gamma(1-\theta)}
\quad \mbox{as }\ n\to+\infty, $$
to obtain 
\begin{multline*}
0\le\int_\D \Big[{}_2F_1(1-\theta,-\theta;1;|w|^2)-
{}_2F_1(1-\theta,-\theta;1;1)\Big]\,|g(w)|^2\,\diff A_{\alpha-2\theta}(w)
\\ = 
\int_\D \sum_{n=1}^{+\infty} 
\frac{|(-\theta)_n|(1-\theta)_n}{(n!)^2}(1-|w|^{2n})
\,|g(w)|^2\,\diff A_{\alpha-2\theta}(w) \\
\le\theta\int_\D\sum_{n=1}^{+\infty}\frac{\big[(1-\theta)_n\big]^2}
{(n-\theta)(n!)^2}n^\theta\,
(1-|w|^2)^\theta\,|g(w)|^2\,\diff A_{\alpha-2\theta}(w) 
\\ \le  
C_2(\alpha,\theta)\left(\sum_{n=1}^{+\infty}
\frac{n^{-\theta}}{n-\theta}\right)
\|g\|^2_{\alpha-\theta},
\end{multline*}
for some appropriate positive constant $C_2(\alpha,\theta)$. By putting
$$C_1(\alpha,\theta)=C_2(\alpha,\theta)\sum_{n=1}^{+\infty}
\frac{n^{-\theta}}{n-\theta},$$
the assertion follows, at least for $0<\theta<1$. The remaining case 
$\theta=1$ is trivial. 
\end{proof}

We remark that the assertion of Lemma \ref{lm-7} remains valid if 
on the right hand side of the estimate we replace the squared norm 
$\|g\|^2_{\alpha-\theta}$ by $\|g\|^2_{\alpha-2\theta+\nu}$, for a fixed
number $\nu$ in the interval $0<\nu<2\theta$. 

It follows from Lemma \ref{lm-7} that (\ref{main-4}) can be written as
\begin{equation}
\label{main-6}
\|L_\theta(z,w)g(w)\|^2_{\alpha,-\theta}= 
\frac{\alpha+1}{\theta(\alpha-2\theta+1)}
\left[1-\frac{\Gamma(2\theta+1)}{2[\Gamma(\theta+1)]^2}\right]
\|g\|^2_{\alpha-2\theta} + 
O\big(\|g\|^2_{\alpha-\theta}\big),
\end{equation}
where the constant in the big ``Oh'' term only depends on $\alpha$ and 
$\theta$. 
To proceed in our calculation of the norm of
$$\big\|\Pop^\perp_{\alpha,-\theta}
\big[L_\theta(z,w)\,g(w)\big]\big\|^2_{\alpha,-\theta},$$
we should like to know the norm of the analytic projection of
the function $L_\theta(z,w)\,g(w)$. We do this by calculating the norm of
of each contribution in the expansion of the function around the diagonal,
in accordance with (\ref{normdecomp}) and Proposition \ref{prop-7.11}.

\begin{prop} 
For $g\in\calH_{\alpha-2\theta}(\D)$, we have 
$$\oslash\left[\frac{\Pop_{\alpha,-\theta;N}\big[L_\theta(z,w)\,g(w)\big]}
{(z-w)^N}\right](z)
=\frac{(-1)^{N+1}(1-\theta)_{N+1}}{(N+1)!\,(\alpha+N+2-2\theta)_{N+1}}
\,\,g^{(N+1)}(z),\qquad z\in\D.$$
\label{prop-7}
\end{prop}

\begin{proof}
In view of (\ref{eq-8}),
\begin{multline*}
\oslash \left[\frac{\Pop_{\alpha,-\theta;N}\big[L_\theta(z,w)\,g(w)\big]}
{(z-w)^N}\right](z)=
\thpar(\alpha,-\theta+N)\\
\times\int_\D\int_\D
\frac{(\bar z'-\bar w')^N}
{(1-z\bar z')^{-\theta+N+2}(1-z\bar w')^{\alpha-\theta+N+2}}
\,L_\theta(z',w')\,g(w')\,
|z'-w'|^{-2\theta}\diff A(z')\,\diff A_\alpha(w').
\end{multline*}
We first integrate with respect to $z'$, that is, we compute
\begin{equation*}
\int_\D\frac{(\bar z'-\bar w')^N}
{(1-z\bar z')^{-\theta+N+2}(1-z\bar w')^{\alpha-\theta+N+2}}
\,L_\theta(z',w')\,|z'-w'|^{-2\theta}\diff A(z').
\end{equation*} 
The change of variables
$$z'=\frac{w'+\zeta}{1+\bar w'\zeta}, \qquad \zeta
=\frac{w'-z'}{1-\bar w' z'},$$
leads to
\begin{multline*}
\int_\D\frac{(\bar z'-\bar w')^N}
{(1-z\bar z')^{-\theta+N+2}(1-z\bar w')^{\alpha-\theta+N+2}}
\,L_\theta(z',w')\,|z'-w'|^{-2\theta}\diff A(z')
\\
=\bar w'\frac{(1-|w'|^2)^{N+1-2\theta}}{(1-z\bar w')^{\alpha-2\theta+2N+4}}
\int_\D\frac1{\bar w'\zeta}\,\Big[\big(1+\bar w'\zeta\big)^{\theta-1}-1\Big]
\left(1+\bar\zeta\,\frac{w'-z}{1-z\bar w'}\right)^{\theta-N-2}
\bar\zeta^N\,\frac{\diff A(\zeta)}{|\zeta|^{2\theta}}\\
=\bar w'\,\frac{(1-|w'|^2)^{N+1-2\theta}}{(1-z\bar w')^{\alpha-2\theta+2N+4}}
\sum_{n=0}^{+\infty}
\frac{\left(\begin{array}{c}\theta-1\\N+n+1\end{array}\right)
\left(\begin{array}{c}\theta-N-2\\n\end{array}\right)}{N+n+1-\theta}\,
(\bar w')^{N+n}\,\left(\frac{w'-z}{1-z\bar w'}\right)^n.
\end{multline*}
The integration  with respect to $w'$ then gives
\begin{multline*}
\int_\D\int_\D\frac{(\bar z'-\bar w')^N}
{(1-z\bar z')^{-\theta+N+2}(1-z\bar w')^{\alpha-\theta+N+2}}
\,L_\theta(z',w')\,g(w')\,|z'-w'|^{-2\theta}\diff A(z')\,\diff A_\alpha(w')
\\
=(\alpha+1)\sum_{n=0}^{+\infty}
\frac{\left(\begin{array}{c}\theta-1\\N+n+1\end{array}\right)
\left(\begin{array}{c}\theta-N-2\\n\end{array}\right)}{N+n+1-\theta}\\
\times\int_\D
(\bar w')^{N+n+1}\,(w'-z)^n\,g(w')\,
\frac{(1-|w'|^2)^{N+1+\alpha-2\theta}}{(1-z\bar w')^{\alpha-2\theta+2N+n+4}}
\,\diff A(w').
\end{multline*}
Next, we notice that by differentiating the reproducing identity for
the weighted Bergman kernel $k$ times, we obtain
\begin{equation*}
\int_\D \frac{(\bar w')^k}{(1-z\,\bar w')^{\gamma+k+2}}\,f(w')\,
\diff A_\gamma(w')=
\frac{1}{(\gamma+2)_k}\,f^{(k)}(z);
\end{equation*}
as we implement this into the above identity, the result is
\begin{multline*}
\oslash\left[\frac{\Pop_{\alpha,-\theta;N}\big[L_\theta(z,w)\,g(w)\big]}
{(z-w)^N}\right](z)
\\
=(-1)^{N+1}\,(\alpha+1)\,\frac{\thpar(\alpha,-\theta+N)}{(N+1)!}\,g^{(N+1)}(z)
\sum_{n=0}^{+\infty}\frac{(1-\theta)_{N+n+1}(N+1-\theta)_{n}}{(N+1-\theta)
\,n!\,(\alpha+N-2\theta+2)_{N+n+2}},
\end{multline*}
so that 
\begin{multline*}
\oslash\left[\frac{\Pop_{\alpha,-\theta;N}\big[L_\theta(z,w)\,g(w)\big]}
{(z-w)^N}\right](z)\\
=(-1)^{N+1}\,(\alpha+1)\,\frac{\thpar(\alpha,-\theta+N)}{(N+1)!}\,g^{(N+1)}(z)
\frac{(1-\theta)_N}{(\alpha+N-2\theta+2)_{N+2}}\\ \times
\sum_{n=0}^{+\infty}\frac{(N+1-\theta)_{n}(N+2-\theta)_{n}}
{n!\,(\alpha+2N-2\theta+4)_{n}}\\
=(-1)^{N+1}\,(\alpha+1)\,\frac{\thpar(\alpha,-\theta+N)}{(N+1)!}\,g^{(N+1)}(z)
\frac{(1-\theta)_N}{(\alpha+N-2\theta+2)_{N+2}}\\ \times
\,{}_2F_1\big(N+1-\theta,N+2-\theta;\alpha+2N-2\theta+4;1\big).
\end{multline*}
 If we use (\ref{Gausshypergeom})
as well as Lemma \ref{lm-4}, the proof is completed.
\end{proof}

\begin{cor}
For $g\in\calH_{\alpha-2\theta}(\D)$, we have 
\begin{multline}
\label{main-7}
\Big\|\Pop_{\alpha,-\theta}\big[L_\theta(z,w)\,g(w)\big]
\Big\|^2_{\alpha,-\theta}\\=\sum_{N=0}^{+\infty}
\frac1{\thpar(\alpha,-\theta+N)}\,
\left[\frac{(1-\theta)_{N+1}}{(N+1)!\,(\alpha+N+2-2\theta)_{N+1}}\right]^2
\,\big\|g^{(N+1)}\big\|^2_{\alpha-2\theta+2N+2},
\end{multline}
where the constant $\thpar(\alpha,-\theta+N)$ is as in Lemma \ref{lm-4}.
\label{cor-9.99}
\end{cor}

The next proposition is crucial for our further analysis. 

\begin{prop}
\label{main-asympt} 
$(-1<\alpha<+\infty)$ Fix the real parameter $\nu$, with $0<\nu\le1$.
Then there exists a positive constant $C_3(\alpha,\nu)$ such that 
for each function $g\in\calH_\alpha(\D)$ and every integer 
$n=1,2,3,\ldots$, 
\begin{equation*}
0\le (\alpha+2)_{2n}\|g\|^2_{\alpha}-\big\|g^{(n)}\big\|^2_{\alpha+2n} \le
C_3(\alpha,\nu)\,n^{2\nu}(\alpha+2)_{2n}\|g\|^2_{\alpha+\nu}. 
\end{equation*}
\label{prop-mainasymptidentity}
\end{prop}

\begin{proof}
The first step is to note that the norm in $\calH_\alpha(\D)$ can be 
expressed as follows in terms of the Taylor coefficients: 
$$\|g\|^2_\alpha= \sum_{k=0}^{+\infty}
\frac{k!}{(\alpha+2)_k}\,|\widehat g(k)|^2. $$
We then have
\begin{multline*}
(\alpha+2)_{2n}\|g\|^2_{\alpha}-\big\|g^{(n)}\big\|^2_{\alpha+2n} \\ 
= (\alpha+2)_{2n}\sum_{k=0}^{+\infty}
\frac{k!}{(\alpha+2)_k}\,|\widehat g(k)|^2 - 
\sum_{k=n}^{+\infty}\frac{(k-n)!}{(\alpha+2+2n)_{k-n}}\,\big[(k-n+1)_n\big]^2
|\widehat g(k)|^2 = \\ = 
(\alpha+2)_{2n}\sum_{k=0}^{+\infty}
\left(1-\frac{(k-n+1)_n}{(k+\alpha+2)_n}\right)
\frac{k!}{(\alpha+2)_k}|\widehat g(k)|^2. 
\end{multline*}
The assertion of the proposition follows from this identity together with 
the following technical inequality: 
\begin{equation}
\label{main-9}
0\le 1-\frac{(k-n+1)_n}{(k+\alpha+2)_n}\le C_4(\alpha)\frac{n^{2\nu}}
{(k+1)^\nu},\qquad k=0,1,2,3,\ldots,\,\,\, n=1,2,3,\ldots. 
\end{equation}
The left hand side of this inequality is obvious. The right hand side 
is also more or less obvious (with $C_4(\alpha)=1$) for $k\le  n^2-1$ . 
So, we assume that $k\ge n^2$. Then we have, by the standard properties of 
the logarithm function,
\begin{multline*}
1-\frac{(k-n+1)_n}{(k+\alpha+2)_n}\le
\log\left[\frac{(k+\alpha+2)_n}{(k-n+1)_n} \right] = \\ =  
\sum_{l=1}^n\bigg[\log\left(1+\frac{\alpha+1+l}k\right) - 
\log\left(1-\frac{n-l}k\right)\bigg]\le 
\sum_{l=1}^n\left[\frac{\alpha+1+l}k+ C_5\frac{n-l}k\right]\le \\ \le 
C_4(\alpha)\,\frac{n^2}{k+1}\le C_4(\alpha)\left[\frac{n^2}{k+1}\right]^\nu,
\end{multline*}
for appropriate values of the positive constants $C_4(\alpha)$ and $C_5$. 
We are done.
\end{proof}

We are now allowed to replace $\|g^{(N+1)}\|^2_{\alpha-2\theta+2N+2}$ 
in each term of (\ref{main-7}) by the expression
$$(\alpha-2\theta+2)_{2N+2}\,\|g\|^2_{\alpha-2\theta},$$
while estimating the remainder as prescribed by Proposition 
\ref{prop-mainasymptidentity} with $\nu=\theta$. In fact, we get convergence
for the estimate of the remainder term so long as $0<\nu<2\theta$.  
After some algebraic manipulations, we then arrive at 
\begin{equation}
\label{main-10}
\Big\|\Pop_{\alpha,-\theta}\big[L_\theta(z,w)\,g(w)\big]
\Big\|^2_{\alpha,-\theta}= 
\fopar(\alpha,\theta)\|g\|^2_{\alpha-2\theta}+ 
O\big(\|g\|^2_{\alpha-\theta}\big), 
\end{equation}
where
\begin{multline}
\label{main-11}
\fopar(\alpha,\theta)=
\frac{(1-\theta)\Gamma(\alpha+2)\Gamma(\alpha+2-2\theta)}
{\Gamma(\alpha+2-\theta)\Gamma(\alpha+3-\theta)}\\
\times
\sum_{N=0}^{+\infty}\big(\alpha+3-2\theta+2N\big)\,
\frac{(1-\theta)_N(2-\theta)_N\big[(\alpha+2-2\theta)_N\big]^2}
{(\alpha+2-\theta)_N(\alpha+3-\theta)_N\big[(N+1)!\big]^2}.
\end{multline}
The series which comes from summing the estimates for the
remainders converges, by the standard asymptotics of the Pochhammer symbol. 
   
The constant $\fopar(\alpha,\theta)$ can be expressed in terms of the 
generalized hypergeometric function ${}_4F_{3}$. We recall its definition: 
$${}_4 F_3\bigg(\begin{array}{cccc}a_1&a_2&a_3&a_4\\{}&b_1&b_2&b_3
\end{array}\bigg|x\bigg)=1+\sum_{n=1}^{+\infty}
\frac{(a_1)_n(a_2)_n(a_3)_n(a_4)_n}{(b_1)_n(b_2)_n(b_3)_n\,n!}\,x^n,$$
wherever the series converges.
By splitting the last factor in the right hand side of (\ref{main-11}) 
as the sum $\alpha+3-2\theta+2N=(\alpha+2-2\theta+N)+(N+1)$, we obtain 
\begin{multline}
\label{main-12}
\fopar(\alpha,\theta)=\frac{(1-\theta)\Gamma(\alpha+2)
\Gamma(\alpha+2-2\theta)}
{\Gamma(\alpha+2-\theta)\Gamma(\alpha+3-\theta)}
\left\{
\frac{(\alpha+1-\theta)(\alpha+2-\theta)}{\theta(1-\theta)
(\alpha+1-2\theta)} - \right. \\ - 
\frac{(\alpha+1-\theta)(\alpha+2-\theta)}{\theta(1-\theta)
(\alpha+1-2\theta)}
{}_4 F_3\bigg(\begin{array}{cccc}-\theta&1-\theta&
\alpha-2\theta+1&\alpha-2\theta+2\\{}&1&\alpha-\theta+1&\alpha-\theta+2
\end{array}\bigg|1\bigg) \\ 
+ \left. 
{}_4 F_3\bigg(\begin{array}{cccc}1-\theta&2-\theta&
\alpha-2\theta+2&\alpha-2\theta+2\\{}&2&\alpha-\theta+2&\alpha-\theta+3
\end{array}\bigg|1\bigg)
\right\}.
\end{multline}
We combine (\ref{main-3}), (\ref{main-6}), and (\ref{main-10}), to
obtain the following expression for the right hand side of (\ref{main-2}): 
$$\left\{\frac{(\alpha+1)\Gamma(2\theta+1)}
{2\theta\,(\alpha-2\theta+1)\big[\Gamma(\theta+1)\big]^2}
+ \fopar(\alpha,\theta)\right\}\,\|g\|^2_{\alpha-2\theta} + 
O\big(\|g\|^2_{\alpha-\theta}\big). $$
On the other hand, the left hand side of (\ref{main-2}) may be likewise
decomposed into a series by the use of Corollary \ref{cor-1.10} and 
Proposition \ref{prop-7}. For $k=0,1,2,3,\ldots$, we introduce  the
analytic functions $\Phi_{k,\theta}$ by the formula
$$\Phi_{k,\theta}(z)=\oslash\big[\bpartial_z^k\Phi_\theta\big](z),\qquad
z\in\D.$$
We arrive at the following statement.

\begin{prop}
$(-1+2\theta<\alpha<+\infty)$ For $g\in\calH_{\alpha-2\theta}(\D)$, we 
have 
\begin{multline*}
\Big\|\Phi_\theta(z,w)\,g(w)+
\Pop_{\alpha,-\theta}\big[L_\theta(z,w)\,g(w)\big]
\Big\|^2_{\alpha,-\theta}\\
=\sum_{N=0}^{+\infty}\frac1{\thpar(\alpha,-\theta+N)}\,
\left\|b_N\,g^{(N+1)}(z)+
\sum_{k=0}^N a_{k,N}\,\bpartial_z^{N-k}\big[\Phi_{k,\theta}(z)\,g(z)\big]
\right\|^2_{\alpha-2\theta+2N+2},
\end{multline*}
where the constant $\thpar(\alpha,-\theta+N)$ is as in Lemma \ref{lm-4}, and
the other constants are given by 
\begin{equation}
\label{main-13}
b_N=\frac{(-1)^{N+1}(1-\theta)_{N+1}}{(N+1)!\,(\alpha-2\theta+N+2)_{N+1}}
\end{equation}
and
\begin{equation}
\label{main-14}
a_{k,N}=\frac{(-1)^{N-k}}{k!(N-k)!}\,\,\,\frac{(-\theta+k+2)_{N-k}}
{(\alpha-2\theta+N+k+3)_{N-k}}.
\end{equation}
\label{prop-6.1}
\end{prop}

Finally, we express the main inequality (\ref{main-2}) in the following 
guise.

\begin{thm}
\label{maininequalitythm}
$(-1+2\theta<\alpha<+\infty)$
There exists a constant $C_6(\alpha,\theta)$ depending only on 
$\theta,\alpha$, with $0<\theta\le1$, such that for any 
$g\in \calH_\alpha(\D)$,
\begin{multline*}
\sum_{N=0}^{+\infty}\frac1{\thpar(\alpha,-\theta+N)}\,
\left\|b_N\,g^{(N+1)}(z)+
\sum_{k=0}^N a_{k,N}\,\bpartial_z^{N-k}\big[\Phi_{k,\theta}(z)\,g(z)\big]
\right\|^2_{\alpha-2\theta+2N+2}\le  \\ \le 
\left[\frac{(\alpha+1)\Gamma(2\theta+1)}
{2\theta(\alpha-2\theta+1)\big[\Gamma(\theta+1)\big]^2}
+ \fopar(\alpha,\theta)\right]\|g\|^2_{\alpha-2\theta}+ 
C_6(\alpha,\theta)\,\|g\|^2_{\alpha-\theta},  
\end{multline*}
where the constants $\thpar(\alpha,N-\theta)$, $ b_N$, 
$a_{k,N}$, and $\fopar(\alpha,\theta)$ are given by 
Lemma \ref{lm-4} and equations $\rm(\ref{main-13})$, $\rm(\ref{main-14})$, 
and {\rm (\ref{main-12})}, respectively. 
\end{thm} 

\section{The algebra of $\varphi$-forms}\label{sec-alg}

In the classical theory of univalent functions, we frequently encounter
expressions like
$$\frac{\varphi''(z)}{\varphi'(z)}\qquad\text{and}\qquad
\frac{\varphi'''(z)}{\varphi'(z)}-
\frac32\left[\frac{\varphi''(z)}{\varphi'(z)}\right]^2,$$
where the first is known as the {\sl logarithmic derivative of the 
derivative} (or the {\sl pre-Schwarz\-ian derivative}), and the second is 
known as the {\sl Schwarzian derivative} of the given univalent function 
$\varphi\in\schlicht$. There are higher-order expressions of a similar 
nature, and it seems reasonable to try to classify them.

An expression of the form
$$\frac{\varphi^{(n+1)}(z)}{\varphi'(z)},$$
with $n$ a positive integer, is said to be a {\sl monomial $\varphi$-form 
of degree $n$ and bidegree $1$}. The degree and bidegree are additive under 
multiplication, which means that, for instance,
$$\frac{\varphi'''(z)\varphi''(z)}{[\varphi'(z)]^2}$$
is a monomial $\varphi$-form of degree $3$ and bidegree $2$. We form 
linear combinations of $\varphi$-forms of the same degree $n$ and the same 
bidegree $k$, and say that the resulting expression is a monomial  
$\varphi$-form of degree $n$ and bidegree $k$.
We may also form linear combinations of monomial $\varphi$-forms of 
the same degree $n$ but of different bidegrees, and speak of the result as a
$\varphi$-form of degree $n$ (without a bidegree). As we form sums of
monomial $\varphi$-forms of various degrees, the maximum of which is $n$,
we get a {\sl $\varphi$-form} with the degree $n$. This way, we get an
{\sl algebra of $\varphi$-forms}. As far as we are concerned, only monomial
$\varphi$-forms will be of any interest.  
\medskip

\noindent{\bf Explicit calculation of the functions $\Phi_{k,\theta}$.} 
We recall the formula
$$\Phi_\theta(z,w)=\frac{1}{z-w}\,\left\{
\frac{\varphi'(z)}{\varphi'(w)}\,
\left(\frac{\varphi(z)-\varphi(w)}{\varphi'(w)\,(z-w)}
\right)^{-\theta-1}-1\right\},\qquad (z,w)\in\D^2,\,\,\,z\neq w.$$
We expand $\varphi(z)$ in a Taylor series about $z=w$:
$$\varphi(z)=\varphi(w)+\sum_{j=1}^{+\infty}\frac{\varphi^{(j)(w)}}{j!}
\,(z-w)^j.$$
This means that
$$\frac{\varphi(z)-\varphi(w)}{\varphi'(w)\,(z-w)}=
\sum_{j=1}^{+\infty}\frac1{j!}\,{\varphi^{(j)(w)}}{\varphi'(w)}
\,(z-w)^{j-1}=1+
\sum_{j=2}^{+\infty}\frac1{j!}\,\frac{\varphi^{(j)}(w)}{\varphi'(w)}
\,(z-w)^{j-1},$$
which leads to 
\begin{multline*}
\left(\frac{\varphi(z)-\varphi(w)}{\varphi'(w)\,(z-w)}\right)^{-\theta-1}
=\bigg[1+
\sum_{j=2}^{+\infty}\frac1{j!}\,\frac{\varphi^{(j)}(w)}{\varphi'(w)}
\,(z-w)^{j-1}
\bigg]^{-\theta-1}\\
=\sum_{n=0}^{+\infty}\left(\begin{array}{c}-\theta-1\\n\end{array}\right)
\,\left(\sum_{j=2}^{+\infty}\frac1{j!}\,
\frac{\varphi^{(j)}(w)}{\varphi'(w)}\,(z-w)^{j-1}\right)^n\\
=1+\sum_{n=1}^{+\infty}\left(\begin{array}{c}-\theta-1\\n\end{array}\right)
\,(z-w)^n\,\left(\sum_{j=2}^{+\infty}\frac1{j!}\,
\frac{\varphi^{(j)}(w)}{\varphi'(w)}\,(z-w)^{j-2}\right)^n.
\end{multline*}
We also have the Taylor series expansion for $\varphi'$, which leads to
$$\frac{\varphi'(z)}{\varphi'(w)}=1+\sum_{k=2}^{+\infty}\frac1{(k-1)!}\,
\frac{\varphi^{(k)}(w)}{\varphi'(w)}\,(z-w)^{k-1}.$$
As we multiply these expressions together, we obtain
\begin{multline*}
\frac{\varphi'(z)}{\varphi'(w)}\,
\left(\frac{\varphi(z)-\varphi(w)}{\varphi'(w)\,(z-w)}
\right)^{-\theta-1}=1+\sum_{k=2}^{+\infty}\frac1{(k-1)!}\,
\frac{\varphi^{(k)}(w)}{\varphi'(w)}\,(z-w)^{k-1}\\
+\sum_{k=1}^{+\infty}\frac1{(k-1)!}\,
\frac{\varphi^{(k)}(w)}{\varphi'(w)}\,(z-w)^{k-1}\\
\times \sum_{n=1}^{+\infty}
\left(\begin{array}{c}-\theta-1\\n\end{array}\right)
\,(z-w)^n\,\left(\sum_{j=2}^{+\infty}\frac1{j!}\,
\frac{\varphi^{(j)}(w)}{\varphi'(w)}\,(z-w)^{j-2}\right)^n,
\end{multline*}
so that
\begin{multline*}
\Phi_\theta(z,w)=\sum_{k=2}^{+\infty}\frac1{(k-1)!}\,
\frac{\varphi^{(k)}(w)}{\varphi'(w)}\,(z-w)^{k-2}\\
+\sum_{k=1}^{+\infty}\frac1{(k-1)!}\,
\frac{\varphi^{(k)}(w)}{\varphi'(w)}\,(z-w)^{k-1}\\
\times \sum_{n=1}^{+\infty}
\left(\begin{array}{c}-\theta-1\\n\end{array}\right)
\,(z-w)^{n-1}\,\left(\sum_{j=2}^{+\infty}\frac1{j!}\,
\frac{\varphi^{(j)}(w)}{\varphi'(w)}\,(z-w)^{j-2}\right)^n.
\end{multline*}
The next step is to note that
\begin{multline*}
\left(\sum_{j=2}^{+\infty}\frac1{j!}\,
\frac{\varphi^{(j)}(w)}{\varphi'(w)}\,(z-w)^{j-2}\right)^n\\=
\sum_{j_1,\ldots,j_n=1}^{+\infty}\frac{\varphi^{(j_1+1)}(w)\cdots
\varphi^{(j_n+1)}(w)}{(j_1+1)!\cdots (j_n+1)!\,\,
[\varphi'(w)]^n}\,(z-w)^{j_1+\ldots+j_n-n},
\end{multline*}
so that we get
\begin{multline}
\Phi_\theta(z,w)=\sum_{l=0}^{+\infty}\frac1{(l+1)!}\,
\frac{\varphi^{(l+2)}(w)}{\varphi'(w)}\,(z-w)^{l}\\
+\left\{1+\sum_{l=1}^{+\infty}\frac1{l!}\,
\frac{\varphi^{(l+1)}(w)}{\varphi'(w)}\,(z-w)^{l}\right\}\\
\times \sum_{n=1}^{+\infty}
\left(\begin{array}{c}-\theta-1\\n\end{array}\right)
\sum_{j_1,\ldots,j_n=1}^{+\infty}\frac{\varphi^{(j_1+1)}(w)\cdots
\varphi^{(j_n+1)}(w)}{(j_1+1)!\cdots (j_n+1)!\,\,
[\varphi'(w)]^n}\,(z-w)^{j_1+\ldots+j_n-1}.
\label{eq-5.01}
\end{multline}

For integers $k,n$, with $1\le n\le k$, we introduce the function
$$\Psi_{k,n}(z)=\sum_{(j_1,\ldots,j_n)\in I(k,n)}
\frac{\varphi^{(j_1+1)}(z)\cdots
\varphi^{(j_n+1)}(z)}{(j_1+1)!\cdots (j_n+1)!\,\,[\varphi'(z)]^n},$$
where $I(k,n)$ is the set of all $n$-tuples $(j_1,\ldots,j_n)$ of positive
integers with $j_1+\ldots+j_n=k$. We realize that $\Psi_{k,n}(z)$ is 
a monomial $\varphi$-form of degree $k$ and bidegree $n$. We calculate that,
for instance,
$$\Psi_{k,1}(z)=\frac{\varphi^{(k+1)}(z)}{(k+1)!\,\varphi'(z)},\qquad
\Psi_{k,2}(z)=\sum_{l=1}^{k-1}
\frac{\varphi^{(l+1)}(z)\,\varphi^{(k-l+1)}(z)}{(l+1)!(k-l+1)!
\,[\varphi'(z)]^2}.$$

\begin{prop} For $k=0,1,2,\ldots$, we have
\begin{equation*}
\Phi_{k,\theta}(z)=
\oslash\big[\bpartial_z^k\Phi_\theta\big](z)
=(k+1-\theta)\,k!\sum_{n=1}^{k+1}
\frac{(-1)^{n-1}(\theta+1)_{n-1}}{n!}\Psi_{k+1,n}(z).
\end{equation*}
\label{prop-5}
\end{prop}

\begin{proof}
We calculate that
\begin{multline*}
\sum_{l=1}^{+\infty}\frac1{l!}\,
\frac{\varphi^{(l+1)}(w)}{\varphi'(w)}\,(z-w)^{l}\\
\times \sum_{n=1}^{+\infty}
\left(\begin{array}{c}-\theta-1\\n\end{array}\right)
\sum_{j_1,\ldots,j_n=1}^{+\infty}\frac{\varphi^{(j_1+1)}(w)\cdots
\varphi^{(j_n+1)}(w)}{(j_1+1)!\cdots (j_n+1)!\,\,
[\varphi'(w)]^n}\,(z-w)^{j_1+\ldots+j_n}\\
=\sum_{n=1}^{+\infty}
\left(\begin{array}{c}-\theta-1\\n\end{array}\right)
\sum_{j_0,j_1,\ldots,j_n=1}^{+\infty}(j_0+1)\,\frac{\varphi^{(j_0+1)}(w)\cdots
\varphi^{(j_n+1)}(w)}{(j_0+1)!\cdots (j_n+1)!\,\,
[\varphi'(w)]^n}\,(z-w)^{j_0+\ldots+j_n},
\end{multline*}
and realize that the expression involving the sum over $j_0,\ldots,j_n$ is
essentially of the same type as the sum appearing on the previous line which
was over $j_1,\ldots,j_n$. By (\ref{eq-5.01}), then, the $k$-th order Taylor
coefficient is 
\begin{multline*}
\frac1{k!}\,\Phi_{k,\theta}(w)=\frac1{k!}\,
\oslash\big[\bpartial_z^k\Phi_\theta\big](w)
=\frac1{(k+1)!}\,
\frac{\varphi^{(k+2)}(w)}{\varphi'(w)}\\
+\sum_{n=1}^{k+1}
\left(\begin{array}{c}-\theta-1\\n\end{array}\right)
\sum_{(j_1,\ldots,j_n)\in I(k+1,n)}\frac{\varphi^{(j_1+1)}(w)\cdots
\varphi^{(j_n+1)}(w)}{(j_1+1)!\cdots (j_n+1)!\,\,
[\varphi'(w)]^n}\\
+\sum_{n=1}^{k}\left(\begin{array}{c}-\theta-1\\n\end{array}\right)
\sum_{(j_0,\ldots,j_n)\in I(k+1,n+1)}(j_0+1)\,\frac{\varphi^{(j_0+1)}(w)
\cdots\varphi^{(j_n+1)}(w)}{(j_0+1)!\cdots (j_n+1)!\,\,
[\varphi'(w)]^n}.
\end{multline*}
We see that
\begin{multline*}
\sum_{(j_0,\ldots,j_n)\in I(k+1,n+1)}(j_0+1)\,\frac{\varphi^{(j_0+1)}(w)
\cdots\varphi^{(j_n+1)}(w)}{(j_0+1)!\cdots (j_n+1)!\,\,
[\varphi'(w)]^n}\\
=\frac{n+k+2}{n+1}\sum_{(j_0,\ldots,j_n)\in I(k+1,n+1)}
\frac{\varphi^{(j_0+1)}(w)
\cdots\varphi^{(j_n+1)}(w)}{(j_0+1)!\cdots (j_n+1)!\,\,
[\varphi'(w)]^n}\\
=\frac{n+k+2}{n+1}\,\Psi_{k+1,n+1}(w),
\end{multline*}
which leads to the simplification
\begin{multline*}
\frac1{k!}\,\Phi_{k,\theta}(w)=\frac1{(k+1)!}\,
\frac{\varphi^{(k+2)}(w)}{\varphi'(w)}
+\sum_{n=1}^{k+1}
\left(\begin{array}{c}-\theta-1\\n\end{array}\right)\,\Psi_{k+1,n}(w)\\
+\sum_{n=1}^{k}\left(\begin{array}{c}-\theta-1\\n\end{array}\right)\,
\frac{n+k+2}{n+1}\,\Psi_{k+1,n+1}(w).
\end{multline*}
As we change the order of summation a bit, and change variables from $w$ to
$z$, the assertion of the proposition follows.
\end{proof}

\begin{remark} \rm
It follows that the expression $\Phi_{k,\theta}(z)$ is a monomial 
$\varphi$-form of degree $k+1$. 
\end{remark}
\medskip

\noindent{\bf Derivatives of powers of $\varphi'$.}
Let $\lambda$ be a complex parameter, and consider the function
$$g_\lambda(z)=[\varphi'(z)]^\lambda=\exp\big[\lambda\log\varphi'(z)\big],
\qquad z\in\D,$$
where $\log\varphi'(z)$ takes the value $0$ at $z=0$, and is analytic
throughout the disk $\D$. We compute that
\begin{equation}
g'_\lambda(z)=\lambda\,\frac{\varphi''(z)}{\varphi'(z)}\,g_\lambda(z),
\label{eq-7.1}
\end{equation}
and
\begin{equation}
g''_\lambda(z)=\lambda\left(\frac{\varphi'''(z)}{\varphi'(z)}
+(\lambda-1)\left[\frac{\varphi''(z)}{\varphi'(z)}\right]^2
\right)\,g_\lambda(z).
\label{eq-7.2}
\end{equation}
Let $\formderiv_{k,\lambda}(z)$ be the function defined by
\begin{equation}
g^{(k)}_\lambda(z)=\formderiv_{k,\lambda}(z)\,g_\lambda(z),
\label{eq-7.3}
\end{equation}
which means that
$$\formderiv_{1,\lambda}(z)=\lambda\,\frac{\varphi''(z)}{\varphi'(z)},
\qquad \formderiv_{2,\lambda}(z)=\lambda\frac{\varphi'''(z)}{\varphi'(z)}
+\lambda(\lambda-1)\left[\frac{\varphi''(z)}{\varphi'(z)}\right]^2.$$
From the rules of differentiation, we have that
$$\formderiv_{k+1,\lambda}(z)=\formderiv_{k,\lambda}'(z)+\lambda\,
\frac{\varphi''(z)}{\varphi'(z)}\,\formderiv_{k,\lambda}(z).$$
This allows us to successively calculate a few higher order factors 
$\formderiv_{k,\lambda}(z)$, such as $\formderiv_{3,\lambda}(z)$:
\begin{equation}
\formderiv_{3,\lambda}(z)=\lambda\frac{\varphi^{(4)}(z)}{\varphi'(z)}
+3\lambda(\lambda-1)\,\frac{\varphi'''(z)\varphi''(z)}{[\varphi'(z)]^2}
+\lambda(\lambda-1)(\lambda-2)
\left[\frac{\varphi''(z)}{\varphi'(z)}\right]^3,
\label{eq-7.4}
\end{equation}
To obtain the formula for the general case, we use the tentative 
representation
\begin{equation}
\formderiv_{k,\lambda}(z)=\sum_{n=1}^{k}(\lambda-n+1)_n
\sum_{(j_1,\ldots,j_n)\in I(k,n)}c(j_1,\ldots,j_n)\,
\frac{\varphi^{(j_1+1)}(z)\cdots\varphi^{(j_n+1)}(z)}
{[\varphi'(z)]^n},
\label{eq-7.5}
\end{equation}
where as before, $I(k,n)$ is the set of all $n$-tuples $(j_1,\ldots,j_n)$ 
of positive integers with $j_1+\ldots+j_n=k$. Also, we assume that the
as of yet undetermined coefficients $c(j_1,\ldots,j_n)$ are invariant under 
permutations, so that, for instance, $c(j_1,\ldots,j_n)=c(j_n,\ldots,j_1)$.
Let ${\mathfrak P}(j_1,\ldots,j_n)$ denote the collection of all 
(different) permutations of the given $n$-tuple $(j_1,\ldots,j_n)$. 
We begin by setting $c(1)=1$, and we define
$$c(j_1,\ldots,j_{n-1},0)=\frac1n\,c(j_1,\ldots,j_{n-1}),$$
for positive integers $j_1,\ldots,j_{n-1}$. All the other values of the 
constants appearing in (\ref{eq-7.5}) are obtained iteratively from the
formula
$$c(j_1,\ldots,j_n)=\frac{n}{\big|{\mathfrak P}(j_1,\ldots,j_n)\big|}
\sum_{(J_1,\ldots,J_n)\in{\mathfrak P}(j_1,\ldots,j_n)} c(J_1,\ldots,J_{n-1},
J_n-1),$$ 
where the absolute value sign is used to denote the number of elements. 

\begin{remark} \rm
For all $k=1,2,3,\ldots$, the expression $\formderiv_{k,\lambda}$ is a 
monomial $\varphi$-form of degree $k$.
\end{remark}

\section{Estimates of the integral means spectrum}

\noindent{\bf An estimate based on the first diagonal term.} 
In this section, we shall use the first term on the left hand side of
the inequality of Theorem \ref{maininequalitythm} to obtain an estimate of
the universal integral means spectrum $B_\schlicht(\tau)$, which is of
interest mainly for $\tau\in\C$ near the origin. 

Throughout this section, we assume that $\varphi$ is a sufficiently smooth
function of the class $\schlicht$; to make this precise, we shall suppose
that $\varphi$ is analytic and univalent in slightly larger disk than $\D$. 
For appropriate values of the real parameter $\beta$ 
($\beta$ is allowed to depend on $\tau$), we shall obtain estimates of 
the norms $\big\|(\varphi')^{\tau/2}\big\|_{\beta-1}$ that
are  uniform in $\varphi$. By a standard dilation argument, we then get the
same uniform norm estimate for general $\varphi\in\schlicht$ as well. In view
of (\ref{eq:formula-spectrum}), this leads to the estimate
$B_\schlicht(\tau)\le \beta$. 

The following proposition is based on Theorem \ref{maininequalitythm}, with
only the first term on the left hand side counted. It uses a fixed value
for the parameter $\theta$. For the formulation, we need the expression
\begin{equation}
\fipar(\beta,\theta)=\frac{(\beta+2\theta)\,\Gamma(2\theta+1)}
{2\,\theta \beta\,[\Gamma(\theta+1)]^2}
+ \fopar(\beta+2\theta -1, \theta),
\label{eq:K-bound}
\end{equation}
where the function $\fopar$ is as in $(\ref{main-11})$ or 
$(\ref{main-12})$.

\begin{prop}
\label{means-1term} Fix $\tau\in\C\setminus\{0\}$ and $\theta$ with 
$0<\theta<1$. Suppose that for some positive real $\beta$, the following 
inequality holds: 
\begin{equation}
\fipar(\beta,\theta)
<(1-\theta)(\beta+1)(\beta+2)\left|\frac 1{\beta+1}-\frac 1\tau \right|^2
\frac{\Gamma(\beta+1+2\theta)\,\Gamma(\beta+2)}
{\Gamma(\beta+1+\theta)\,\Gamma(\beta+2+\theta)},
\label{condition-1term}
\end{equation}
where the function $\fipar$ is as above. Suppose, in addition, that 
$$\big\|(\varphi')^{\tau/2}\big\|_{\beta-1+\theta}=O(1)$$ 
holds uniformly in $\varphi\in\schlicht$. Then we also have 
$$\big\|(\varphi')^{\tau/2}\big\|_{\beta-1}=O(1)$$
uniformly in $\varphi\in\schlicht$. In particular, 
${\mathrm B}_\schlicht(\tau)\le \beta$. 
\end{prop}

\begin{proof}
If we take into account only the first term of the sum on the left hand 
side of the inequality in Theorem \ref{maininequalitythm}, and pick
$\alpha=\beta+2\theta-1$, we obtain
\begin{multline}
\label{means-1}
\frac 1{\thpar(\beta+2\theta-1,-\theta)}\left\|-\frac{1-\theta}{\beta+1}\,g'+ 
\frac{1-\theta}2\,\frac{\varphi''}{\varphi'}\,g\right\|^2_{\beta+1} \\
\le \,
\fipar(\beta,\theta)\,\|g\|^2_{\beta-1}+ O\Big(\|g\|^2_{\beta-1+\theta}\Big),
\end{multline}
for an arbitrary $g\in\calH_{\beta-1}(\D)$. 
Here, we used the fact that 
$$\Phi_{0,\theta}(z)=\oslash\Phi_\theta(z)={\frac{1-\theta}2 \,
\frac{\varphi''(z)}{\varphi'(z)}},\qquad z\in\D,$$
which is an almost trivial case of Proposition \ref{prop-5}. 

The next step is to apply the estimate (\ref{means-1}) to the functions
$$g(z)=g_\tau(z)=\big[\varphi'(z)\big]^{\tau/2},$$
and to make the observation that
\begin{equation}
\frac{\varphi''(z)}{\varphi'(z)}\,g_\tau(z) = \frac 2\tau \,g'_\tau(z),
\qquad z\in\D.
\label{eq-idder}
\end{equation}
By Proposition \ref{main-asympt} (with $\nu=\theta$), we have 
$$\big\|g'\big\|^2_{\beta+1}=(\beta+1)(\beta+2)\,\|g\|^2_{\beta-1}+ 
O\big(\|g\|^2_{\beta-1+\theta}\big)$$
holds generally, so that if we combine it with the above observation and
recall the formula of Lemma \ref{lm-4}, we obtain from (\ref{means-1}) that
\begin{multline*}
\left\{(1-\theta)
\frac{\Gamma(\beta+1+2\theta)\Gamma(\beta+2)}
{\Gamma(\beta+1+\theta)\Gamma(\beta+2+\theta)}
(\beta+1)(\beta+2)\left|\frac 1{\beta+1}-\frac 1\tau \right|^2- 
\fipar(\beta,\theta)\right\}\,\big\|g_\tau\big\|^2_{\beta-1}\\ = 
O\big(\|g_\tau\|^2_{\beta-1+\theta}\big),
\end{multline*}
which implies the assertion of the proposition. 
\end{proof}

\begin{remark}\rm
A part of the assertion of Proposition \ref{means-1term}, namely 
${\mathrm B}_\schlicht(\tau)\le\beta$, remains true under
the weaker assumption of ``$\le$'' in (\ref{condition-1term}). This is
so because in the case of equality in (\ref{condition-1term}) for given 
$\theta$, $\beta$, and $\tau$, we may move $\tau$ slightly so as to
achieve ``$<$''. Using the continuity of the function 
${\mathrm B}_\schlicht$, the asserted inequality follows by taking the limit.
\label{means-rk0}
\end{remark}

We may use the above proposition iteratively to obtain successively better
bounds for the function $B_\schlicht(\tau)$ starting from some some trivial
bound, like what follows from the pointwise K\oe{}be-Bieberbach 
estimate (\ref{eq-koebeest}). A more general estimate is
\begin{equation}
\Big|\big[\varphi'(z)\big]^\tau\Big|\le
\frac{(1+|z|)^{2|\tau|-\re\tau}}{(1-|z|)^{2|\tau|+\re\tau}},\qquad z\in\D,
\label{eq:esttau}
\end{equation}
which works for general $\tau\in\C$; it is obtained if we integrate 
(\ref{eq:Bieberbach}), to get
$$\Big|\log\varphi'(z)+\log\big(1-|z|^2\big)\Big|\le
2\log\frac{1+|z|}{1-|z|},\qquad z\in\D,$$
and perform the appropriate algebraic manipulations. It follows from
(\ref{eq:esttau}) that for fixed $\tau\in\C$,
\begin{equation}
\Big\|\big[\varphi'\big]^{\tau/2}\Big\|^2_{2|\tau|+\re\tau-1+\eps}
=O(1)  
\label{means-2}
\end{equation}
holds uniformly in $\varphi\in\schlicht$, for all positive values of $\eps$.
\medskip

\noindent{\bf A first estimate of the integral means spectrum near the
origin.}   We apply Proposition \ref{means-1term} to obtain asymptotic 
bounds for the function $B_\schlicht(t)$ for $t$ near the origin. 

\begin{prop}
\label{means-origin1} Fix a $\theta$ with $0<\theta<1$. We then have
\begin{equation}
\label{means-3}
\limsup_{\C\ni\tau\to 0}\,\frac{B_\schlicht(\tau)}{|\tau|^2} \le 
\frac{1+\theta}{2(1-\theta)}.
\end{equation}
\end{prop}

\begin{proof}
Pick a positive $\eps$, and let 
$$\beta=\beta(\tau)=
\left[\frac{1+\theta}{2(1-\theta)}+\eps\right]\,|\tau|^2. $$
We plug this $\beta$ into both sides of (\ref{condition-1term}),
and observe that the left hand side behaves like
$$\frac{2(1-\theta)\Gamma(2\theta+1)}{\Gamma(1+\theta)\Gamma(2+\theta)}\,
\frac 1{|\tau|^2}+o\left(\frac 1{|\tau|^2}\right) \quad\text{ as }\,\,\,
\tau\to 0, $$
while the right hand side behaves like
$$\left[\frac{1+\theta}{2(1-\theta)}+\eps\right]^{-1}
\frac{\Gamma(2\theta+1)}{\left[\Gamma(\theta+1)\right]^2}\,\frac 1{|\tau|^2}
+o\left(\frac 1{|\tau|^2}\right) \quad \text{ as }\,\,\, \tau\to 0,$$
which shows that condition (\ref{condition-1term}) is fulfilled for 
sufficiently small values of $|\tau|$.  As the trivial estimate
(\ref{means-2}) show that
$$\Big\|\big[\varphi'\big]^{\tau/2}\Big\|^2_{\beta(\tau)-1+\theta}=O(1)$$ 
for sufficiently small $|\tau|$, we may apply Proposition \ref{means-1term} 
to deduce that 
$$\Big\|\big[\varphi'\big]^{\tau/2}\Big\|^2_{\beta(\tau)-1}=O(1)$$
holds uniformly in $\varphi$ for sufficiently small $|\tau|$. The desired 
assertion follows. 
\end{proof}

\begin{cor}
\label{cor:means-origin2}
We have 
\begin{equation*}
\limsup_{\C\ni\tau\to 0}\,\frac{B_\schlicht(\tau)}{|\tau|^2} \le \frac 12. 
\end{equation*}
\end{cor}

\begin{proof}
Let $\theta\to 0^+$ in (\ref{means-3}).
\end{proof}
\medskip

\noindent{\bf The improved estimate of the integral means spectrum near the
origin.}   Below, we  obtain a better constant instead of $\frac12$ in the 
estimate of Corollary \ref{cor:means-origin2}. 

Naturally, if we take into account more terms of the sum in the
left hand side of the inequality in Theorem \ref{maininequalitythm}, 
we obtain more precise information. We now analyze the estimate
obtained by considering the first two terms. As in the proof of 
Proposition \ref{means-1term}, we fix some $\theta$ with $0<\theta<1$ and 
some positive $\beta$, and we plug in $\alpha=\beta+2\theta-1$ and 
$g=g_\tau=[\varphi']^{\tau/2}$ into Theorem \ref{maininequalitythm}, throwing
away all but the first two terms on the left hand side. We use 
Proposition \ref{prop-5} to evaluate $\Phi_{k,\theta}(z)$ for $k=0,1$, and
the identity (\ref{eq-idder}) to obtain,
for $0<\beta<+\infty$, 
\begin{multline}
\label{means-2termsineq}
(1-\theta)(\beta+1)(\beta+2)\,
\frac{\Gamma(\beta+1+2\theta)\Gamma(\beta+2)}
{\Gamma(\beta+1+\theta)\Gamma(\beta+2+\theta)}
\,\left|\frac 1{\beta+1}-
\frac 1\tau\right|^2\Big\|\big[\varphi'\big]^{\tau/2}\Big\|^2_{\beta-1}\\ + 
(2-\theta)\,\frac{\Gamma(\beta+2\theta+1)\Gamma(\beta+4)}
{\Gamma(\beta+\theta+2)\Gamma(\beta+\theta+3)}\times 
%\mbox{\huge {\bf X}}
\clubsuit
\\
\le 
\fipar(\beta,\theta)\,\Big\|\big[\varphi'\big]^{\tau/2}\Big\|^2_{\beta-1}
+O\Big(\Big\|\big[\varphi'\big]^{\tau/2}\Big\|^2_{\beta+\theta-1}\Big),
\end{multline}
where
\begin{multline}
\label{means-2termsineqX}
%\mbox{\huge {\bf X}}
\clubsuit
= \left\|\frac{1-\theta }{2(\beta+2)(\beta+3)}\, \dz^2
\Big\{\big[\varphi'\big]^{\tau/2}\Big\}- 
\frac{1-\theta}{2(\beta+3)}\,
\dz\left\{\frac{\varphi''}{\varphi'}\,\big[\varphi'\big]^{\tau/2}\right\} 
\right. \\ + \left.   
\left\{\frac 16\,\frac{\varphi'''}{\varphi'}-\frac{\theta+1}8
\left(\frac{\varphi''}{\varphi'}\right)^2\right\}
\big[\varphi'\big]^{\tau/2} \right\|^2_{\beta+3},
\end{multline}
and $\dz=\diff/{\diff z}$ stands for the operator of differentiation. 
As before, we first apply this inequality to estimate $B_\schlicht(\tau)$ 
near the origin. We consider $\beta=\beta(\tau)=B_0\,|\tau|^2$, where 
$B_0$ is some fixed  constant with $0<B_0<\frac12$. 
We put $\theta=\theta(\tau)=4|\tau|$, and plug these values into
(\ref{means-2termsineq}). By the trivial estimate (\ref{means-2}), we
have 
$$\Big\|\big[\varphi'\big]^{\tau/2}\Big\|^2_{\beta(\tau)+\theta(\tau)-1}
=O(1),$$
uniformly in $\varphi\in\schlicht$ for each fixed $\tau\in\C$. 
Then (\ref{means-2termsineq}) takes the following form:
\begin{multline}
\label{means-4}
\frac{2+\epsilon_1(\tau)}{|\tau|^2}\,
\Big\|\big[\varphi'\big]^{\tau/2}\Big\|^2_{-1+\beta(\tau)} \\
 + \big(6+\epsilon_2(\tau)\big)\left\|\bigg\{
\left(\frac 1{24}+\epsilon_3(\tau)\right)
\left[\frac{\varphi''}{\varphi'}\right]^2 + 
\epsilon_4(\tau)\,\frac{\varphi'''}{\varphi'}\bigg\}
\big[\varphi'\big]^{\tau/2}\right\|^2_{3+\beta(\tau)} 
\\ 
\le\frac{1+\epsilon_5(\tau)}{B_0\,|\tau|^2}\,\Big\|\big[\varphi'\big]^{\tau/2}
\Big\|^2_{-1+\beta(\tau)} + O(1),
\end{multline}
where the last $O(1)$ is uniform in $\varphi\in\schlicht$ for each fixed 
$\tau$. For $k=1,2,3,4,5$, the functions $\epsilon_k(\tau)$ satisfy
$$\lim_{\tau\to 0}\epsilon_k(\tau)=0;$$
and for $k=1,2,5$, the functions are in addition real-valued.

\begin{lemma}
\label{means-phi'''}
 $(-1<\alpha<+\infty)$
There exists a positive constant $C_7(\alpha)$ such that, for any 
$g\in\calH_\alpha(\D)$,
$$\left\|\frac{\varphi'''}{\varphi'}\,g\right\|^2_{\alpha+4}\le 
C_7(\alpha)\,\|g\|^2_\alpha. $$
Moreover, 
$$C_7(\alpha)=O\left(\frac1{\alpha +1}\right)$$ 
as $\alpha\to-1^+$. 
\end{lemma}

\begin{proof} 
The assertion follows from the identity 
$$\frac{\varphi'''}{\varphi'}\,g= \frac \diff{\diff z}
\left(\frac{\varphi''}{\varphi'}\,g\right)
+\left[\frac{\varphi''}{\varphi'}\right]^2g -\frac{\varphi''}{\varphi'}\,g'$$
combined with the classical pointwise estimate (\ref{eq:Bieberbach})
and Proposition \ref{main-asympt}. 
\end{proof}

As we apply the above lemma, we obtain from (\ref{means-4}) that 
the inequality 
\begin{equation}
\label{means-5}
\bigg\|\left[\frac{\varphi''}{\varphi'}
\right]^2\big[\varphi'\big]^{\tau/2}\bigg\|^2_{\beta(\tau)+3}
\le 
\frac{96}{|\tau|^2}\,\left(\frac 1{B_0}-2+\eps\right)
\Big\|\big[\varphi'\big]^{\tau/2}\Big\|^2_{\beta(\tau)-1} + O(1)
\end{equation}
holds for each fixed positive $\eps$, for sufficiently small values of
$|\tau|$. 
 
\begin{lemma}
\label{means-CB}
$(0<\beta<+\infty)$ For each $g\in\calH_{\beta-1}$, we have 
\begin{equation}
\label{means-CBineq}
\left\|\frac{\varphi''}{\varphi'}\,g\right\|^2_{\beta+1}\le 
\frac{\beta+2}{\sqrt{\beta(\beta+4)}}\,\|g\|_{\beta-1} \,
\bigg\|\left[\frac{\varphi''}{\varphi'}\right]^2g\bigg\|_{\beta+3}. 
\end{equation}
\end{lemma}

\begin{proof}
This follows from a standard application of the 
Cauchy-Schwarz--Bunyakovski\u\i{} inequality. 
\end{proof}

By estimate (\ref{means-5}) and Lemma \ref{means-CB},
we have the following chain of inequalities (as before, 
$\beta(\tau)=B_0|\tau|^2$): 
\begin{multline*}
\Big\|\big[\varphi'\big]^{\tau/2}\Big\|^2_{\beta(\tau)-1}= 
\frac {|\tau|^2}{4\,(\beta(\tau)+1)(\beta(\tau)+2)}\,
\left\|\frac{\varphi''}{\varphi'}\,\big[\varphi'\big]^{\tau/2}
\right\|^2_{1+\beta(\tau)} + O(1) 
\\
\le \frac {|\tau|^2}{4\,(\beta(\tau)+1)\sqrt{\beta(\tau)(\beta(\tau)+4)}}\, 
\bigg\|\left[\frac{\varphi''}{\varphi'}\right]^2 
\big[\varphi'\big]^{\tau/2}\bigg\|_{\beta+3}\,
\Big\|\big[\varphi'\big]^{\tau/2}\Big\|_{\beta(\tau)-1} + O(1) 
\\
\le \frac {1+\epsilon_6(\tau)}{8\sqrt{B_0}}\,
\sqrt{96\left(\frac 1{B_0}-2+\eps\right)}\,
\Big\|\big[\varphi'\big]^{\tau/2}\Big\|^2_{\beta(\tau)-1}
+ O\bigg(\Big\|\big[\varphi'\big]^{\tau/2}\Big\|_{\beta(\tau)-1}\bigg)+O(1),
\end{multline*}
where the function $\epsilon_6(\tau)$ is real-valued with limit
$\epsilon_6(\tau)\to0$ as $\tau\to0$. 
This inequality implies that
$$\Big\|\big[\varphi'\big]^{\tau/2}\Big\|_{\beta(\tau)-1}=O(1)$$ 
uniformly in $\varphi\in\schlicht$, provided that 
$$\frac{1+\epsilon_6(\tau)}{8\sqrt{B_0}}\,
\sqrt{96\left(\frac 1{B_0}-2+\eps\right)}<1. $$
We conclude that 
$$\limsup_{\C\ni \tau\to 0}\,\frac{B_\schlicht(\tau)}{|\tau|^2}\le B_0$$
holds for each real constant $B_0$, $0<B_0<\frac12$, for which
$$96\,\left(\frac 1{B_0}-2\right)<64\,B_0. $$
By solving this last inequality for $B_0$, we obtain the following estimate.

\begin{thm}
\label{means-origin}
We have that
\begin{equation} 
\label{means-originineq}
\limsup_{\C\ni\tau\to 0}\frac{B_\schlicht(\tau)}{|\tau|^2} \le
\frac{\sqrt{15}-3}{2} = 0.43649\ldots 
\end{equation}
\end{thm}

\begin{remark} \rm The best previous estimate of this type was 
$B_\schlicht (t)\le (3+\eps)\, t^2$ for real $t$ near the origin 
(see \cite{Pomm}). 
\end{remark}
\medskip

\noindent{\bf An optimization method to estimate ${\mathrm B}_\schlicht$ 
using two terms.}
Our next goal is to estimate the function $B_\schlicht(\tau)$ using the
the inequality (\ref{means-2termsineq}), which employs the first two terms
on the left hand side of the inequality in Theorem \ref{maininequalitythm}.
This time we intend to take into account somehow all possible values of 
$\theta$ at the same time, rather than considering a single value at a time. 
This of course requires that the estimates we have obtained so far 
are sufficiently uniform in $\theta$, if $\theta$ is confined to some
compact interval $[\theta_0,1]$, which is true and possible to verify
without too much effort. 
We fix $\tau\in\C$ and $\beta$ with $0<\beta<+\infty$, and
rewrite (\ref{means-2termsineq}) as follows, using (\ref{eq-7.1}): 
\begin{multline}
\label{means-2termsopt}
\left\| A_1(\theta)\,\dz^2\Big\{\big[\varphi'\big]^{\tau/2}\Big\} + 
A_2(\theta)\left[\frac{\varphi''}{\varphi'}\right]^2
\big[\varphi'\big]^{\tau/2} \right\|_{\beta+3} \\ 
\le \,\Big\|\big[\varphi'\big]^{\tau/2}\Big\|_{-1+\beta} 
+ O\Big(\Big\|\big[\varphi'\big]^{\tau/2}\Big\|_{-1+\beta+\theta}\Big),
\end{multline}
where
\begin{multline}
\label{means-A1}
A_1(\theta)=\left[\frac{1-\theta}{2(\beta+2)(\beta+3)} - 
\frac{1-\theta}{\tau\,(\beta+3)} + \frac 1{3\tau} \right]
\bigg\{(2-\theta)\,\frac{\Gamma(\beta+2\theta+1)\,\Gamma(\beta+4)}
{\Gamma(\beta+\theta+2)\,\Gamma(\beta+\theta+3)}\bigg\}^{\frac12} \\
\times\left\{\fipar(\beta,\theta)- (1-\theta)(\beta+1)(\beta+2)
\frac{\Gamma(\beta+2\theta+1)\,\Gamma(\beta+2)}
{\Gamma(\beta+\theta+1)\,\Gamma(\beta+\theta+2)}
\left|\frac 1{\beta+1}-\frac 1\tau\right|^2
\right\}^{-1/2},
\end{multline}
and
\begin{multline}
\label{means-A2}
A_2(\theta)=\left[\frac 16\left(1-\frac \tau2\right)-\frac{\theta+1}8
\right]
\left\{(2-\theta)\,\frac{\Gamma(\beta+2\theta+1)\Gamma(\beta+4)}
{\Gamma(\beta+\theta+2)\Gamma(\beta+\theta+3)}\right\}^{\frac12}
\\
\times\left\{\fipar(\beta,\theta)- (1-\theta)(\beta+1)(\beta+2)
\frac{\Gamma(\beta+2\theta+1)\Gamma(\beta+2)}
{\Gamma(\beta+\theta+1)\Gamma(\beta+\theta+2)}
\left|\frac 1{\beta+1}-\frac 1\tau\right|^2\right\}^{-\frac12};
\end{multline}
we recall the definition of the function $\fipar(\beta,\theta)$ in 
(\ref{eq:K-bound}). Without loss of generality, we may assume that 
\begin{equation}
\label{means-lastcondition}
(1-\theta)(\beta+1)(\beta+2)
\frac{\Gamma(\beta+2\theta+1)\Gamma(\beta+2)}
{\Gamma(\beta+\theta+1)\Gamma(\beta+\theta+2)}
\left|\frac 1{\beta+1}-\frac 1\tau\right|^2<\fipar(\beta,\theta)
\end{equation}
holds for all $\theta$, $0<\theta\le1$; for otherwise, we may
apply Proposition \ref{means-1term} in conjunction with Remark 
\ref{means-rk0} to get the desired inequality 
${\mathrm B}_\schlicht(\tau)\le \beta$. This means that the square roots
which are used to define the functions $A_1$ and $A_2$ produce real-valued
functions on the whole interval $0<\theta\le1$. 
For each $\theta$, $0<\theta\le1$, we consider the disk 
\begin{equation}
\label{means-Itheta}
\calD_\theta=\bigg\{\,w\in\C\,:\,\,
\big|A_1(\theta)-w\,A_2(\theta)\big| \le 
\frac1{\sqrt{(\beta+1)_4}}\bigg\}. 
\end{equation} 
Here, of course, $(\beta+1)_4=(\beta+1)(\beta+2)(\beta+3)(\beta+4)$.
\medskip

We have the following result.

\begin{prop}
\label{means-optim}
Suppose that there exists a certain $\theta_0$, with $0<\theta_0\le1$, such
that
\smallskip
 
\noindent $(a)$ 
the intersection $\bigcap_{\theta_0<\theta\le1}\calD_\theta$ 
is empty, and 
\smallskip

\noindent $(b)$ the estimate 
$\big\|[\varphi']^{\tau/2}\big\|_{-1+\beta+\theta_0}=O(1)$ holds uniformly in 
$\varphi\in\schlicht$. 
\smallskip

\noindent Then 
$$\big\|[\varphi']^{\tau/2}\big\|_{-1+\beta}=O(1)$$ 
holds uniformly in $\varphi\in\schlicht$, so that in particular,   
${\mathrm B}_\schlicht (\tau)\le\beta$. 
\end{prop}

\begin{proof}
A standard compactness argument shows that the assumption $(a)$ 
remains valid if we replace the disks $\calD_\theta$ by the slightly 
bigger disks 
$$\calD^\eps_\theta = \bigg\{\,w\in\C\,:\,\,
\big|A_1(\theta)-w\,A_2(\theta)\big| \le 
\frac{1+\eps}{\sqrt{(\beta+1)_4}}\bigg\},$$
for a small enough positive $\eps$. This means that
$$\inf_{w\in\C}\,\big\| A_1-w\,A_2\big\|_{C[\theta_0,1]}
\ge \frac {1+\eps}{\sqrt{(\beta+1)_4}}$$
holds, if, as is standard, $C[\theta_0,1]$ is the Banach space of 
complex-valued functions continuous in $[\theta_0,1]$, supplied with 
the uniform norm. 
By standard duality, this entails that there exists a complex Borel measure 
$\mu$ on the interval $[\theta_0,1]$ such that the total variation of $\mu$ 
is $1$, and, in addition, 
$$\frac {1+\eps}{\sqrt{(\beta+1)_4}}
\le\left|\int_{\theta_0}^1 A_1(\theta)\,\diff\mu(\theta)\right|
\qquad\text{while}\qquad
\int_{\theta_0}^1 A_2(\theta)\,\diff\mu(\theta)=0. $$
We find that an application of (\ref{means-2termsopt}) leads to
\begin{multline*}
\frac {1+\eps}{\sqrt{(\beta+1)_4}}
\left\|\dz^2\left\{\big[\varphi'\big]^{\tau/2}\right\}
\right\|_{\beta+3}\le
\left\| \left\{\int_{\theta_0}^1 A_1(\theta)\,\diff\mu(\theta)\right\}
\dz^2\left\{\big[\varphi'\big]^{\tau/2}\right\}\right\|_{\beta+3}
\\ = 
\left\|\int_{\theta_0}^1 \left\{
A_1(\theta)\,
\dz^2\left\{\big[\varphi'\big]^{\tau/2}\right\}
+A_2(\theta)\left[\frac{\varphi''}{\varphi'}\right]^2
(\varphi')^{\tau/2} \right\}\, \diff\mu(\theta) \right\|_{\beta+3} \\
\le
\int_{\theta_0}^1 \left\|
A_1(\theta)\,
\dz^2\left\{\big[\varphi'\big]^{\tau/2}\right\}
+A_2(\theta)\left[\frac{\varphi''}{\varphi'}\right]^2
(\varphi')^{\tau/2} \,\right\|_{\beta+3} |\diff\mu(\theta)|
\\
\le \, 
\Big\|\big[\varphi'\big]^{\tau/2}\Big\|_{-1+\beta}+ 
O\Big(\Big\|\big[\varphi'\big]^{\tau/2}\Big\|_{-1+\beta+\theta_0}\Big). 
\end{multline*}
In view of Proposition \ref{main-asympt} and the assumption $(b)$, the 
desired conclusion follows.
\end{proof}

The moral content of Proposition \ref{means-optim} is that {\sl we are able 
to obtain the estimate 
$$\big\|[\varphi']^{\tau/2}\big\|_{-1+\beta}=O(1)$$ 
uniformly over all $\varphi\in\schlicht$ for as long as the criterion
\begin{equation}
\bigcap_{0<\theta\le1}\calD_\theta(\secpar,\tau)=\emptyset
\label{intersec=empty}
\end{equation}
is fulfilled}, where $\calD_\theta(\secpar,\tau)=\calD_\theta$ is as in
(\ref{means-Itheta}). In a concrete situation, of course, we have to start
with a trivial {\sl a priori} estimate, and inch our way down in the scale of
$\beta$'s in accordance with details specified by Proposition 
\ref{means-optim}. We should mention that by Helly's intersection theorem,
(\ref{intersec=empty}) holds if and only if 
\begin{equation*}
\calD_{\theta_1}(\secpar,\tau)\cap\calD_{\theta_2}(\secpar,\tau)\cap
\calD_{\theta_3}(\secpar,\tau) =\emptyset
\end{equation*} 
for some triplet $\theta_1,\theta_2,\theta_3$ with $0<\theta_j\le1$,
$j=1,2,3$.

\begin{remark} \rm For $\tau=t$ real, it suffices to verify the
assumption $(a)$ of Proposition \ref{means-optim} along the real line only,
as can be seen from the observation that the functions $A_1(\theta)$ and
$A_2(\theta)$ are real-valued then. This means that if we put
\begin{equation}
\label{means-Itheta'}
\calI_\theta=\bigg\{\,x\in\R\,:\,\,
\big|A_1(\theta)-x\,A_2(\theta)\big| \le 
\frac1{\sqrt{(\beta+1)_4}}\bigg\}, 
\end{equation} 
which constitutes a closed interval, it is enough to check that
\begin{equation}
\bigcap_{\theta_0\le\theta\le1}\calI_\theta=\emptyset.
\label{means-inters}
\end{equation}
This criterion can be easily checked by computer calculations. 
Indeed, if we denote the left and right end points of $\calI_\theta$ by 
$\alpha_1(\theta)$ and $\alpha_2(\theta)$, so that
$$\calI_\theta=[\alpha_1(\theta),\alpha_2(\theta)],$$
then the criterion (\ref{means-inters}) is equivalent to 
$$\inf_{\theta\in[\theta_0,1]}\alpha_2(\theta)<
\sup_{\theta\in[\theta_0,1]}\alpha_1(\theta), $$
which is easily treated numerically. 
\label{means-rk}
\end{remark}

\begin{remark} \rm
It would be desirable to change the implementation of
the optimization method so that we may incorporate the information
supplied by Lemma \ref{means-CB}, so as to obtain a more optimal estimate
based on the first two terms. If we do this in a straightforward manner,
focussing on the term containing $A_2(\theta)$ instead of $A_1(\theta)$,
we are to replace the intervals $\calI_\theta=\calD_\theta\cap\R$ by
$$\calJ_\theta=\bigg\{\,x\in\R\,:\,\,
\big|A_2(\theta)-x\,A_1(\theta)\big| \le 
\frac{t^2}{4(\beta+1)\sqrt{\beta(\beta+4)}}\bigg\},$$
and the criterion 
$\bigcap_{\theta_0\le\theta\le1}\calJ_\theta=\emptyset$
then permits us to conclude that ${\mathrm B}_\schlicht(t)\le\beta$.
Numerical simulation shows that this criterion is more powerful for
(real) $t$ near the origin than the criterion $(a)$ of Proposition 
\ref{means-optim}. 
\label{remark-cb}
\end{remark}
\medskip

\noindent{\bf Numerical implementation.} By successive application of 
Proposition \ref{means-optim} for real $\tau=t$, taking into account
Remark \ref{means-rk}, we obtain the estimate ${\mathrm B}_\schlicht(t)\le
{\mathrm B}_*(t)$, where the function ${\mathrm B}_*(t)$ is tabulated 
below. We use suitably small values of $\theta_0$. The function 
${\mathrm B}_*(t)$ is also graphed. For some values of $t$, the method
outlined in Remark \ref{remark-cb} is used in place of Proposition 
\ref{means-optim}; this is then indicated with an asterisk (*).

The tabulated bounds for ${\mathrm B}_\schlicht(-1)$ and 
${\mathrm B}_\schlicht(-2)$ are to be compared with the bounds that were 
found recently by the second-named author in \cite{Shi}; there, it was 
shown that ${\mathrm B}_\schlicht(-1)\le0.420$ and 
${\mathrm B}_\schlicht(-2)\le1.246$. It should be noted that the inequality 
of Theorem 1 in \cite{Shi} leading to these bounds is a particular case 
of our main inequality -- the inequality of Theorem \ref{maininequalitythm} 
-- if we put $\theta=1$ and, like in (\ref{means-2termsopt}), take into
account  only the first two terms in the sum on the left hand side. In 
this particular case, the first term vanishes and the constant 
$C_6(\alpha,\theta)$  which appears in (\ref{maininequalitythm}) vanishes
as well, because $L_\theta=0$ for $\theta=1$. 
\medskip

\begin{center}
\begin{tabular}{|c c  c|}\hline
 $t$ & ${\mathrm B}_*(t)$   &  $\max\{-t-1,0\}$       \\ \hline
    $-20.000$     &  $19.028$     & $19.000$       \\ \hline
    $-10.000$      &  $9.040$     &  $9.000$       \\ \hline
    $-8.000$      &  $7.049$      &  $7.000$       \\ \hline
    $-6.000$      &  $5.067$      &  $5.000$       \\ \hline
    $-5.000$      &  $4.082$      &  $4.000$       \\ \hline
    $-4.000$      &  $3.105$      &  $3.000$       \\ \hline
    $-3.000$      &  $2.144$      &  $2.000$       \\ \hline
    $-2.500$      &  $1.674$      &  $1.500$       \\ \hline
    $-2.400$      &  $1.582$      &  $1.400$       \\ \hline
    $-2.300$      &  $1.490$      &  $1.300$       \\ \hline
    $-2.200$      &  $1.398$      &  $1.200$       \\ \hline
    $-2.100$      &  $1.308$      &  $1.100$       \\ \hline
    $-2.000$      &  $1.218$      &  $1.000$       \\ \hline
    $-1.900$      &  $1.130$      &  $0.900$       \\ \hline
    $-1.800$      &  $1.042$      &  $0.800$       \\ \hline
    $-1.752$      &  $1.001$      &  $0.752$       \\ \hline
    $-1.700$      &  $0.956$      &  $0.700$       \\ \hline
    $-1.600$      &  $0.871$      &  $0.600$       \\ \hline
    $-1.500$      &  $0.787$      &  $0.500$       \\ \hline
    $-1.400$      &  $0.706$      &  $0.400$       \\ \hline
    $-1.300$      &  $0.626$      &  $0.300$       \\ \hline
    $-1.200$      &  $0.549$      &  $0.200$       \\ \hline
    $-1.100$      &  $0.474$      &  $0.100$       \\ \hline
    $-1.000$      &  $0.403$      &  $0.000$       \\ \hline
    $-0.900$      &  $0.336$      &  $0.000$       \\ \hline
    $-0.800$      &  $0.272$      &  $0.000$       \\ \hline
    $-0.700$      &  $0.213^*$    &  $0.000$       \\ \hline
    $-0.600$      &  $0.159^*$    &  $0.000$       \\ \hline
    $-0.500$      &  $0.112^*$    &  $0.000$       \\ \hline
    $-0.400$      &  $0.072^*$    &  $0.000$       \\ \hline
    $-0.300$      &  $0.0404^*$   &  $0.000$       \\ \hline
    $-0.200$      &  $0.0179^*$   &  $0.000$       \\ \hline
    $-0.150$      &  $0.0100^*$   &  $0.000$       \\ \hline
    $-0.100$      &  $0.00443^*$  &  $0.000$       \\ \hline
    $-0.050$      &  $0.00110^*$  &  $0.000$       \\ \hline
\end{tabular}
$\qquad$
\begin{tabular}{|c c c|}\hline
 $t$ & ${\mathrm B}_*(t)$   & $\max\{3t-1,0\}$        \\ \hline
    $0.000$      &  $0.00000$     &  $0.000$          \\ \hline
    $0.050$      &  $0.00141^*$   &  $0.000$        \\ \hline
    $0.100$      &  $0.0065$      &  $0.000$        \\ \hline
    $0.150$      &  $0.0157$      &  $0.000$        \\ \hline
    $0.200$      &  $0.031$     &  $0.000$          \\ \hline
    $0.250$      &  $0.056$     &  $0.000$          \\ \hline
    $0.300$      &  $0.101$     &  $0.000$          \\ \hline
    $0.350$      &  $0.190$     &  $0.050$          \\ \hline
    $0.400$      &  $0.314$     &  $0.200$          \\ \hline
    $0.450$      &  $0.447$     &  $0.350$          \\ \hline
    $0.500$      &  $0.585$     &  $0.500$          \\ \hline
    $0.600$      &  $0.870$     &  $0.800$          \\ \hline
    $0.700$      &  $1.159$     &  $1.100$          \\ \hline
    $0.800$      &  $1.452$     &  $1.400$          \\ \hline
    $0.900$      &  $1.746$     &  $1.700$          \\ \hline
    $1.000$      &  $2.041$     &  $2.000$          \\ \hline
    $1.100$      &  $2.337$     &  $2.300$          \\ \hline
    $1.200$      &  $2.634$     &  $2.600$          \\ \hline
    $1.300$      &  $2.932$     &  $2.900$          \\ \hline
    $1.400$      &  $3.230$     &  $3.200$          \\ \hline
    $1.500$      &  $3.528$     &  $3.500$          \\ \hline
    $1.600$      &  $3.826$     &  $3.800$          \\ \hline
    $1.700$      &  $4.124$     &  $4.100$          \\ \hline
    $1.800$      &  $4.423$     &  $4.400$          \\ \hline
    $1.900$      &  $4.722$     &  $4.700$          \\ \hline
    $2.000$      &  $5.021$     &  $5.000$          \\ \hline
    $2.100$      &  $5.320$     &  $5.300$          \\ \hline
    $2.200$      &  $5.619$     &  $5.600$          \\ \hline
    $2.300$      &  $5.918$     &  $5.900$          \\ \hline
    $2.400$      &  $6.217$     &  $6.200$          \\ \hline
    $2.500$      &  $6.517$     &  $6.500$          \\ \hline
    $3.000$      &  $8.014$     &  $8.000$          \\ \hline
    $4.000$      &  $11.011$    &  $11.000$         \\ \hline
    $5.000$      &  $14.010$     &  $14.000$        \\ \hline
    $6.000$      &  $17.010$     &  $17.000$        \\ \hline
\end{tabular}
\\[0.5cm]
{\bf TABLES 1 AND 2.}
\end{center}

\vfill\eject

\begin{picture}(100,200)(-100,1)
\put(-110,0){\vector(1,0){380}}
\put(75,0){\vector(0,1){200}}
\put(-75,-5){\line(0,1){4}}
\put(-79,-12){$-2$}
\put(0,-5){\line(0,1){4}}
\put(-4,-12){$-1$}
\put(73,-12){$0$}
\put(150,-5){\line(0,1){4}}
\put(148,-12){$1$}
\put(225,-5){\line(0,1){4}}
\put(223,-12){$2$}
\put(75,75){\line(1,0){4}}
\put(83,72){$1$}
\put(75,150){\line(1,0){4}}
\put(83,147){$2$}
\put(275,-12){$t$}
\put(83,190){$B$}
\curve(-100,100, 0,0)
\put(-85,25){$B=-t-1$}
\curve(100,0, 150,150)
\put(140,100){$B=3t-1$}
\put(0,30){$B={\mathrm B}_*(t)$}
\curve(-90,104.85, -75,91.35, -60,78.15, -45,65.325, -30,52.95, -15,41.175, 
0,30.225, 15,20.4, 30,12.075, 45,5.70,  60,1.5,  75,0,  82.5,0.49,  
90,2.325,  93.75,4.185, 96,5.9197, 97.5,7.575, 101.25,14.22, 103.5,19.725,
105,23.55, 112.5,43.875, 120,65.25, 135,108.9, 150,153.075 )
%\curve(75,100, 75,100, 100,100)
%\put(-50,75){\circle*{2.5}}
%\put(0,75){\circle*{2.5}}
%\put(50,75){\circle*{2.5}}
%\put(100,75){\circle*{2.5}}
%\put(150,75){\circle*{2.5}}
%\put(250,115){\circle*{2.5}}
\thicklines
\end{picture}
\bigskip
\medskip

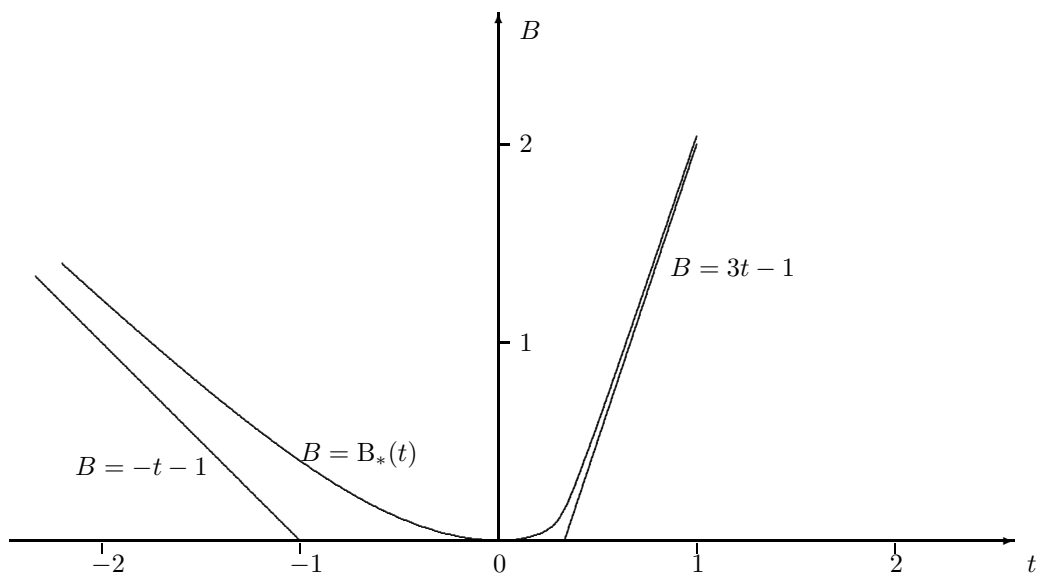
\noindent{{\bf FIGURE 1.} Graph of $B={\mathrm B}_*(t)$, the estimated 
universal spectral function; support lines included.}
\bigskip

\begin{remark} \rm By taking advantage of the fact that the function 
${\mathrm B}_\Sigma(t)$ is convex, with 
${\mathrm B}_\Sigma(t)\le{\mathrm B}_\schlicht(t)$ and 
${\mathrm B}_\Sigma(2)=1$, we derive from a somewhat larger supply of
sample values of the graphed function ${\mathrm B}_*(t)$ that 
${\mathrm B}_\Sigma(1)\le 0.4600$, improving the best earlier known 
estimate, due to Makarov and Pommerenke \cite{MakPom}, which was
${\mathrm B}_\Sigma(1)\le 0.4886$. The value of ${\mathrm B}_\Sigma(1)$
describes the growth of the length of Green lines (the level curves of the
Green function) as they approach the boundary of an arbitrary simply 
connected bounded planar domain. It also determines the rate of decay of
the Laurent series coefficients of functions in the class $\Sigma$ (see
\cite{CarJon}). 
\end{remark}
\medskip

\noindent{\bf The optimization method to estimate ${\mathrm B}_\schlicht$ 
using three or more terms.} How do we implement the optimization method
if we take into account more than two terms on the left hand side of 
the inequality of Theorem \ref{maininequalitythm}? 
We outline here briefly an extension of the method which applies to the case
of three terms. The method may of course be extended to include more than
three terms as well.

For simplicity, we consider real $\tau=t$ only. As we take the first three 
terms on the left hand side of the inequality of 
Theorem \ref{maininequalitythm} into account, putting,
as before, $\alpha=\beta+2\theta-1$ and $g=[\varphi']^{t/2}$,  
we obtain an inequality of the form
\begin{multline}
\label{means-3terms}
\left\| A_1(\theta)\,\dz^2\Big\{\big[\varphi'\big]^{t/2}\Big\} + 
A_2(\theta)\left[\frac{\varphi''}{\varphi'}\right]^2
\big[\varphi'\big]^{t/2} \right\|^2_{\beta+3}+
\frac 1{(\beta+5)(\beta+6)}\\
\times
\left\|A_3(\theta)\,\dz^3\left\{\big[\varphi'\big]^{t/2}\right\} + 
A_4(\theta)\,{\dz}
\left\{\left[\frac{\varphi''}{\varphi'}\right]^2
\big[\varphi'\big]^{t/2}\right\} + 
A_5(\theta)\left[\frac{\varphi''}{\varphi'}\right]^3\big[\varphi'\big]^{t/2} 
\right\|^2_{\beta+5}
\,\le \\ \le\, 
\Big\|\big[\varphi'\big]^{t/2}\Big\|^2_{\beta-1} + 
O\bigg(\Big\|\big[\varphi'\big]^{t/2}\Big\|^2_{\beta-1+\theta}\bigg),
\end{multline}
where the functions $A_1$ and $A_2$ are given by (\ref{means-A1}) and 
(\ref{means-A2}), and the functions $A_3$, $A_4$, $A_5$ are continuous 
on $]0,1]$, and given by certain explicit expressions. As before, we 
assume that condition (\ref{means-lastcondition}) is fulfilled for all 
$\theta$, $0<\theta\le1$. The process of deriving equation 
(\ref{means-3terms}) involves not only Theorem \ref{maininequalitythm},
but also some of the algebraic results of Section \ref{sec-alg}. 
A counterpart to Proposition \ref{means-optim} is the following. 

 \begin{prop}
\label{means-optim3terms}
Let $\calE_\theta$ denote the ellipse in $(x,y)$-plane defined by 
the condition 
$$\big|A_1(\theta)-x\,A_2(\theta)\big|^2
+ \big|A_3(\theta)-x\,A_4(\theta)-y\,A_5(\theta)\big|^2 \le 
\frac 1{(\beta+1)_4}. $$
Suppose that there exists a certain $\theta_0$, with $0<\theta_0\le1$, 
such that 
\smallskip

\noindent$(a)$ the intersection 
$\bigcap_{\theta\in[\theta_0,1]}\calE_\theta$ is empty; 
\smallskip

\noindent$(b)$ $\nphip^2_{-1+\beta+\theta_0}=O(1)$ uniformly in 
$\varphi\in\schlicht$. 
\smallskip

\noindent Then 
$$\Big\|\big[\varphi'\big]^{\tau/2}\Big\|_{-1+\beta}=O(1)$$  
holds uniformly in $\varphi\in\schlicht$ and, in particular,   
${\mathrm B}_\schlicht (t)\le \beta$. 
\end{prop}

\begin{proof}
First, we introduce the operator of integration ${\mathbf I}_0$,
$${\mathbf I}_0 f(z)=\int_0^{z}f(w)\,dw,\qquad z\in\D.$$
Then we apply Proposition \ref{main-asympt} to the second term on 
the left-hand side of (\ref{means-3terms}), which allows us to rewrite 
(\ref{means-3terms}) in the form
\begin{multline}
\label{means-3terms'}
\left\| A_1(\theta)\,\dz^2\Big\{(\varphi')^{t/2}\Big\} + 
A_2(\theta)\left[\frac{\varphi''}{\varphi'}\right]^2
\big[\varphi'\big]^{t/2} \right\|^2_{\beta+3} 
\\ +
\left\|A_3(\theta)\,\dz^2\left\{\big[\varphi'\big]^{t/2}\right\} + 
A_4(\theta)\left[\frac{\varphi''}{\varphi'}\right]^2
\big[\varphi'\big]^{t/2} + 
A_5(\theta)\,{\mathbf I}_0\bigg[\left[\frac{\varphi''}{\varphi'}\right]^3
\big[\varphi'\big]^{t/2}\bigg]
\right\|^2_{\beta+3}
\\ \le\, 
\Big\|\big[\varphi'\big]^{t/2}\Big\|^2_{-1+\beta} 
+ O\bigg(\Big\|\big[\varphi'\big]^{t/2}\Big\|^2_{-1+\beta+\theta}\bigg).
\end{multline}
A standard compactness argument shows that the assumption $(a)$ 
remains valid if the ellipses $\calE_\theta$ are replaced by 
slightly larger ellipses $\calE^\eps_\theta$, defined by 
$$\big|A_1(\theta)-xA_2(\theta)\big|^2
+ \big|A_3(\theta)-xA_4(\theta)-yA_5(\theta)\big|^2 \le 
\frac {1+\eps}{(\beta+1)_4},$$
provided that the positive number $\eps$ is small enough. Moreover, a
similar argument shows that we may assume that a finite intersection 
of $\calE^\eps_\theta$ is empty: 
$$\bigcap_{\theta\in \calF}\calE^\eps_\theta=\emptyset,$$
for some finite subset $\calF$ of the interval $[\theta_0,1]$. This 
condition is equivalent to having
$$\max_{\theta\in \calF}
\left\{\big|A_1(\theta)-xA_2(\theta)\big|^2
+\big|A_3(\theta)-xA_4(\theta)-yA_5(\theta)\big|^2\right\} 
> \frac{1+\eps}{(\beta+1)_4}$$
for all $x,y\in\R$, or, expressed differently,
\begin{equation}
\label{means-slut}
\operatorname{dist}_\calX\left[
\left(\begin{matrix} A_1 \\ A_3\end{matrix} \right),\,\, 
\spann\left\{ 
\left(\begin{matrix} A_2 \\ A_4 \end{matrix} \right),
\left(\begin{matrix} 0  \\ A_5 \end{matrix} \right)
\right\}\right] >
\sqrt{\frac{1+\eps}{(\beta+1)_4}},
\end{equation}
where ``$\spann$'' means the $\R$-linear span, and $\operatorname{dist}_\calX$
is the distance function on the space $\calX$, the $\R$-linear 
space of vector-valued functions 
$$\theta\mapsto \left(\begin{matrix}\xi_1(\theta) \\
\xi_2(\theta)\end{matrix}\right),\qquad \theta\in \calF, $$
supplied with the norm 
$$\left\| 
\left(\begin{matrix}\xi_1 \\
\xi_2\end{matrix}\right)\right\|_\calX = 
\max_{\theta\in \calF} 
\sqrt{|\xi_1(\theta)|^2+|\xi_2(\theta)|^2}.$$
The $\R$-linear space $\calX^*$ of vector-valued functions 
$$\theta\mapsto \left(\begin{matrix}\mu_1(\theta) \\
\mu_2(\theta)\end{matrix}\right), \qquad\theta\in \calF,$$
supplied with the norm 
$$\left\| 
\left(\begin{matrix}\mu_1 \\
\mu_2\end{matrix}\right)\right\|_{\calX^*} = 
\sum_{\theta\in \calF}
\sqrt{|\mu_1(\theta)|^2+|\mu_2(\theta)|^2}, $$
is then dual to $\calX$ with respect to the natural dual pairing
$$\left<
\left(\begin{matrix}\xi_1\\
\xi_2\end{matrix}\right),
\left(\begin{matrix}\mu_1\\
\mu_2\end{matrix}\right) \right> =
\sum_{\theta\in \calF} \Big\{\xi_1(\theta)\mu_1(\theta) + 
\xi_2(\theta)\mu_2(\theta) \Big\}. $$
By standard duality theory, the inequality (\ref{means-slut}) 
means that there exists a vector-valued function 
$$\left(\begin{matrix}\mu_1 \\
\mu_2\end{matrix}\right) \in \calX^* $$ 
which satisfies
\begin{equation}
\left\| \left(\begin{matrix}\mu_1\\
\mu_2\end{matrix}\right) \right\|_{\calX^*} =1;\qquad
\sqrt {\frac {1+\eps}{(\beta+1)_4}}<
\left| \left<
\left(\begin{matrix} A_1 \\ A_2 \end{matrix} \right), 
\left(\begin{matrix}\mu_1 \\
\mu_2 \end{matrix}\right) \right>\right| , 
\label{eq:normcontrol}
\end{equation}
while 
$$\left<
\left(\begin{matrix} A_2 \\ A_4 \end{matrix} \right), 
\left(\begin{matrix}\mu_1 \\
\mu_2 \end{matrix}\right) \right>=
\left<\left(\begin{matrix} 0 \\ A_5 \end{matrix} \right), 
\left(\begin{matrix}\mu_1 \\
\mu_2 \end{matrix}\right) \right>=0. $$
We then have 
\begin{multline*}
\sqrt{\frac {1+\eps}{(\beta+1)_4}}\,
\left\|\dz^2\big[\varphi'\big]^{t/2}\right\|_{\beta+3} <
\left\|\sum_{\theta\in \calF} 
\Big\{\mu_1(\theta)A_1(\theta)+\mu_2(\theta)A_3(\theta)
\Big\}\,\dz^2 \big[\varphi'\big]^{t/2} \right\|_{\beta+3} \\ 
= \left\|\sum_{\theta\in \calF}\left\{\mu_1(\theta)\left[A_1(\theta)\,
\dz^2\Big\{\big[\varphi'\big]^{t/2}\Big\} + 
A_2(\theta)\left[\frac{\varphi''}{\varphi'}\right]^2
\big[\varphi'\big]^{t/2} \right] +\mu_2(\theta)\right.  \right.  
\\
\times\Bigg[
A_3(\theta)\,\dz^2\left[(\varphi')^{t/2}\right]
+ A_4(\theta)\left[\frac{\varphi''}{\varphi'}\right]^2\big[\varphi'\big]^{t/2}
+ A_5(\theta)\, {\mathbf I}_0\bigg[\left[\frac{\varphi''}{\varphi'}\right]^3
\big[\varphi'\big]^{t/2}
\Bigg] \Bigg\}\Bigg\|_{\beta+3} \\
\le 
\sum_{\theta\in \calF}\left\{|\mu_1(\theta)|\,\left\|A_1(\theta)\,
\dz^2\Big\{\big[\varphi'\big]^{t/2}\Big\} + 
A_2(\theta)\left[\frac{\varphi''}{\varphi'}\right]^2
\big[\varphi'\big]^{t/2} \right\|_{\beta+3} +|\mu_2(\theta)|\right.  
\\\times\,\Bigg\|
A_3(\theta)\,\dz^2\Big\{\big[\varphi'\big]^{t/2}\Big]
+ A_4(\theta)\left[\frac{\varphi''}{\varphi'}\right]^2\big[\varphi'\big]^{t/2}
+ A_5(\theta)\, {\mathbf I}_0\bigg[\left[\frac{\varphi''}{\varphi'}\right]^3
\big[\varphi'\big]^{t/2}
\Bigg] \Bigg\|_{\beta+3}\Bigg\}
\\
\le\Big\|\big[\varphi'\big]^{t/2}\Big\|_{\beta-1} + 
O\bigg(\Big\|\big[\varphi'\big]^{t/2}\Big\|_{\beta-1+\theta}\bigg),
\end{multline*}
where in the last step, we appeal to the Minkowski
inequality, as well as to (\ref{eq:normcontrol}) and (\ref{means-3terms'}). 
Since $\eps$ is positive, this completes the proof, in view of
Proposition \ref{main-asympt}.
\end{proof}

\begin{remark} \rm
When running computer tests based on the criterion of Proposition 
\ref{means-optim3terms}, it is useful to know that -- by Helly's intersection
theorem -- it enough to check that
\begin{equation*}
\calE_{\theta_1}\cap\calE_{\theta_2}\cap
\calE_{\theta_3} =\emptyset
\end{equation*} 
for some triplet $\theta_1,\theta_2,\theta_3$ with $\theta_0\le\theta_j\le1$,
$j=1,2,3$. We have not yet carried out these computer runs.
\end{remark}

\section{Ways to extend the method}
\label{extend}
 
Is it possible to improve our method so that it leads to better bounds for 
the integral means? We feel that one way to achieve improvement is to try 
and replace our starting point inequality (\ref{eq-1}) by some other, more 
appropriate estimate. 
\medskip

\noindent{\bf Moving the boundary branch point $\infty$ to $\mu$.} Let 
$\mu\in\C\setminus\varphi(\D)$. Then the function 
$$\varphi_\mu(z)=\frac{\mu\varphi(z)}{\mu-\varphi(z)}$$
is again in $\schlicht$, and replacing $\varphi$ by $\varphi_\mu$ 
in (\ref{eq-1}) leads to   
\begin{multline}
\int_\D\left|\frac{\varphi'(z)}{\varphi'(w)}\,
\left(\frac{\mu-\varphi(w)}
{\mu-\varphi(z)}\right)^{1-\theta}
\left(\frac{\varphi'(w)\,(z-w)}
{\varphi(z)-\varphi(w)}\right)^{\theta+1}
-\left(\frac{1-|w|^2}{1-\bar w z}\right)^{1-\theta}\right|^2
\frac{\diff A(z)}{|z-w|^{2\theta+2}}\\
\le
\frac1{\theta}\,(1-|w|^2)^{-2\theta}.
\label{eq-3}
\end{multline}
We introduce the notation
$$\Phi_{\theta,\mu}(z,w)=\frac{1}{z-w}\,\left\{
\frac{\varphi'(z)}{\varphi'(w)}\,
\left(\frac{\mu-\varphi(w)}{\mu-\varphi(z)}\right)^{1-\theta}
\left(\frac{\varphi'(w)\,(z-w)}
{\varphi(z)-\varphi(w)}\right)^{\theta+1}
-1\right\},$$
so that
$$\Phi_{\theta,\mu}(z,w)=\left(\frac{\mu-\varphi(w)}{\mu-\varphi(z)}
\right)^{1-\theta}\,\Phi_\theta(z,w)+
\frac1{z-w}\,\left\{\left(\frac{\mu-\varphi(w)}{\mu-\varphi(z)}
\right)^{1-\theta}-1\right\}.$$
Note that as $\mu$ tends to $\infty$ from inside the complement of 
$\varphi(\D)$, 
$$\Phi_{\theta,\mu}(z,w)\to \Phi_{\theta}(z,w).$$
Also, $\Phi_{1,\mu}(z,w)\equiv\Phi_1(z,w)$. In terms of $\Phi_{\theta,\mu}$, 
estimate (\ref{eq-3}) can be written as follows.

\begin{thm}
Fix $\theta$, $0<\theta\le1$.
Let $\varphi\in{\schlicht}$ be arbitrary, and suppose 
$\mu\in\C\setminus\varphi(\D)$.
Then, for all $w\in\D$,
\begin{equation*}
\int_\D\Big|\Phi_{\theta,\mu}(z,w)+L_\theta(z,w)\Big|^2
\frac{\diff A(z)}{|z-w|^{2\theta}}
\le
\frac1{\theta}\,(1-|w|^2)^{-2\theta},
%\label{eq-3'}
\end{equation*}
with equality if and only if $\varphi$ is a full mapping.
\label{thm-1.1}
\end{thm}

One way to spread out the effect of the point $\mu$ in Theorem \ref{thm-1.1}
is to integrate both sides of the inequality with respect to a probability
measure in the variable $\mu$, supported on $\C\setminus\varphi(\D)$. A
particularly attractive choice of such a measure would be 
the harmonic measure
for the point at the origin.
  
The diagonal restriction of the function $\Phi_{\theta,\mu}$
equals
$$\Phi_{\theta,\mu}(z,z)=
\frac{1-\theta}2\,\left(\frac{\varphi''(z)}{\varphi'(z)}
+\frac{2\,\varphi'(z)}{\mu-\varphi(z)}\right) $$
Note that the $\mu$-average of this function with respect to the harmonic 
measure for the origin equals
$$\frac{1-\theta}2\,\left\{\frac{\varphi''(z)}{\varphi'(z)}
+\frac2{z}-\frac{2\,\varphi'(z)}{\varphi(z)}\right\}=
\frac{1-\theta}2\,\,\frac{d}{\diff z}
\log\frac{z^2\,\varphi'(z)}{[\varphi(z)]^2}.$$
The expression $z^2\varphi'(z)/(\varphi(z))^2$ is essentially the derivative 
of a function from $\Sigma$, if we use the inversion map to go from $\D$
to $\D_e$. So, averaging with respect to $\mu$ in this way may lead to
interesting properties for the class $\Sigma$. 
\medskip

\noindent{\bf Nehari's extension of Prawitz' theorem.} Another way to 
generalize inequality (\ref{eq-1}) is to start with a more general initial 
inequality than Prawitz' estimate, as given by Theorem \ref{thm-1}. 
A polynomial version of Prawitz' theorem was obtained by Nehari \cite{Neh},
 and it can be reformulated as the inequality (valid for $0<\theta<1$)
$$\int_{\D} \bigg|\sum_{j=1}^N c_j\,\Lambda_\theta(z,\zeta_j)\bigg|^2
\,\frac{\diff A(z)}{|z|^{2\theta}}\le \sum_{j,k=1}^N c_j\,
\bar c_{k}\,K_\theta(\zeta_j,\zeta_{k}),$$
where
$$\Lambda_\theta(z,\zeta)=\left[\frac{\varphi(z)}z\right]^{-\theta}
\left[\frac{\varphi(\zeta)}{\zeta}\right]^\theta
\frac{\varphi'(z)}{\varphi(z)-\varphi(\zeta)} - \frac1{z-\zeta}$$
and 
$$K_\theta(\xi,\zeta)=\sum_{n=0}^{+\infty}
\frac{(\xi\bar\zeta)^n}{n+\theta}. $$
Prawitz' theorem is the special case when $N=1$, $\zeta_1=0$, and  
$c_1=1$. As we shift the branching point from the origin to an arbitrary 
point $w\in \D$ by the same procedure as in Section 2, we obtain the 
inequality
\begin{equation}
\label{extend-1}
\int_{\D} \bigg|\sum_{j=1}^N c_j\,\Lambda_{\theta}(z,\zeta_j;w)\bigg|^2\,
\frac{\diff A(z)}{|z-w|^{2\theta}} \le 
\sum_{j,k=1}^N c_j\,\bar c_{k}\,K_{\theta}(\zeta_j,\zeta_{k};w),
\end{equation}
where 
$$\Lambda_{\theta}(z,\zeta;w)=\left[\frac{\varphi(z)-\varphi(w)}{z-w}
\right]^{-\theta}\left[\frac{\varphi(\zeta)-\varphi(w)}{\zeta-w}
\right]^\theta
\frac{\varphi'(z)}{\varphi(z)-\varphi(\zeta)} - 
\bigg[\frac{1-\bar w\zeta}{1-\bar wz}\bigg]^{1-\theta}
\frac 1{z-\zeta}$$
and 
$$K_{\theta}(\xi,\zeta;w)=
\frac 1{(1-\bar w\xi)^\theta(1-w\bar\zeta)^{\theta}}\,\,
\sum_{n=0}^{+\infty}\frac 1{n+\theta}
\left[\frac{(\xi-w)(\bar\zeta-\bar w)}{(1-\bar w\xi)(1-w\bar\zeta)}\right]^n.
$$
The inequality (\ref{eq-1}) results from (\ref{extend-1}) if we set $N=1$, 
$c_1=1$, and $\zeta_1=w$. It is not clear whether this is the
optimal choice of the parameters for the method.  
\medskip

\noindent{\bf A Prawitz-type theorem with two internal branching points.}
Inequalities in the spirit of (\ref{eq-1}) (or, more generally, 
in the spirit of (\ref{extend-1})), can be obtained for expressions
associated with more than one branching point; after all, (\ref{extend-1})
means that we have placed a single (interior) branching point at $w\in\D$.
Let us focus on the case of two (interior) branching points. We place 
one at the origin and the other at the point $w\in\D$,  and supply both
with ``branching multiplicity'' $\frac12$. The inequality that is the
analogue of (\ref{eq-1}) with $\theta=\frac12$ can be shown to assume the form
\begin{equation}
\label{extend-2}
\int_{\D} \big|\Xi(z,w)+N(z,w)\big|^2\,\frac{\diff A(z)}
{|z(z-w)|} \le M(w),\qquad w\in\D,
\end{equation}
where
$$\Xi(z,w)=\left(\frac{\varphi(w)}w\right)^{1/2}
\left(\frac{\varphi(z)}{z}\right)^{-1/2} 
\frac{\varphi'(z)}{\varphi'(w)}\left(\frac{\varphi'(w)(z-w)}
{\varphi(z)-\varphi(w)}\right)^{3/2},$$
and the functions $N(z,w)$ and $M(w)$ are given by certain explicit
expressions involving elliptic functions. Again, we have equality for all 
full mappings $\varphi\in\schlicht$. The method to derive the 
inequality (\ref{extend-2}) as well as its generalizations to arbitrary
$0<\theta<1$ will be explained in forthcoming papers. 
Unfortunately, the inequality (\ref{extend-2}) -- and similar inequalities
with branching points at the origin and at the point $w$ with 
branching multiplicities  $1-\theta$ and $\theta$ ($0<\theta<1$), 
respectively -- do not seem to yield any new information as regards  
integral means spectra. The analysis of these inequalities leads to 
expressions of just the same type as before, with the only difference
that they are written in terms of the meromorphic function $1/\varphi(z)$,
which becomes an element of the class $\Sigma$ after the change of variables 
$z\mapsto 1/z$.

\bigskip

\noindent Hedenmalm and Shimorin: Department of Mathematics, Royal
Institute of Technology, S--100 44 Stockholm, Sweden. Email:
{\tt haakanh@math.kth.se}, {\tt shimorin@math.kth.se}

\end{document}